\input amstex\documentstyle{amsppt}  
\pagewidth{12.5cm}\pageheight{19cm}\magnification\magstep1
\topmatter
\title From conjugacy classes in the Weyl group to unipotent classes\endtitle
\author G. Lusztig\endauthor
\address{Department of Mathematics, M.I.T., Cambridge, MA 02139}\endaddress
\thanks{Supported in part by the National Science Foundation}\endthanks
\endtopmatter   
\document
\define\Irr{\text{\rm Irr}}
\define\sneq{\subsetneqq}
\define\codim{\text{\rm codim}}

\define\ds{\dot s}

\define\ul{\un l}

\define\uG{\un G}
\define\uuG{\un{\un G}}

\define\uWW{\un\WW}

\define\mpb{\medpagebreak}

\define\hs{\hat s}

\define\dsv{\dashv}

\define\pe{\perp}
\define\si{\sim}

\define\sqc{\sqcup}

\define\qua{\quad}

\define\dx{\dot x}

\define\lb{\linebreak}

\define\op{\oplus}

\define\part{\partial}
\define\em{\emptyset}
\define\imp{\implies}

\define\n{\notin}
\define\iy{\infty}
\define\m{\mapsto}
\define\do{\dots}

\define\lra{\leftrightarrow}

\define\sub{\subset}    

\define\T{\times}
\define\ti{\tilde}
\define\nl{\newline}
\redefine\i{^{-1}}

\define\un{\underline}

\define\ot{\otimes}
\define\bbq{\bar{\QQ}_l}

\define\Ad{\text{\rm Ad}}

\define\End{\text{\rm End}}

\define\tr{\text{\rm tr}}

\define\a{\alpha}
\redefine\b{\beta}
\redefine\c{\chi}
\define\g{\gamma}
\redefine\d{\delta}
\define\e{\epsilon}

\define\io{\iota}
\redefine\o{\omega}
\define\p{\pi}
\define\ph{\phi}
\define\ps{\psi}
\define\r{\rho}
\define\s{\sigma}

\define\th{\theta}
\define\k{\kappa}
\redefine\l{\lambda}

\define\x{\xi}

\redefine\G{\Gamma}

\define\Om{\Omega}

\redefine\L{\Lambda}
\define\Ph{\Phi}

\define\kk{\bold k}

\define\nn{\bold n}

\redefine\ss{\bold s}

\define\CC{\bold C}

\define\FF{\bold F}

\define\NN{\bold N}

\define\QQ{\bold Q}

\define\WW{\bold W}
\define\ZZ{\bold Z}

\define\cb{\Cal B}

\define\cf{\Cal F}

\define\ch{\Cal H}

\define\co{\Cal O}
\define\cp{\Cal P}

\define\car{\Cal R}

\define\cv{\Cal V}

\define\cz{\Cal Z}

\define\fa{\frak a}
\define\fb{\frak b}

\define\fg{\frak g}

\define\fn{\frak n}

\define\fB{\frak B}

\define\fV{\frak V}
\define\fZ{\frak Z}

\define\ts{\ti s}

\define\tA{\ti A}

\define\tE{\ti E}

\define\sha{\sharp}

\define\sps{\supset}

\define\CA{Ca}
\define\DL{DL}
\define\EG{EG}
\define\GE{Ge}
\define\GP{GP}
\define\GH{Ch}
\define\KA{Ka}
\define\KL{KL}
\define\LUE{L\"u}
\define\LU{L1}
\define\OR{L2}
\define\CS{L3}
\define\GF{L4}
\define\HEC{L5}
\define\COR{L6}
\define\MI{Mi}
\define\SH{Sh}
\define\SPA{Sp1}
\define\SPS{Sp2}
\define\SKLC{Sp3}
\define\SKL{Sp4}
\define\ST{St}

\head Introduction\endhead
\subhead 0.1\endsubhead
Let $G$ be a connected reductive algebraic group over an algebraically closed field $\kk$ of characteristic 
$p\ge0$. Let $\uG$ be the set of conjugacy classes in $G$. Let $\uuG$ be the set of unipotent conjugacy 
classes in $G$. Let $\uWW$ be the set of conjugacy classes in
the Weyl group $\WW$ of $G$. In \cite{\KL} a (conjecturally injective) map $\uuG@>>>\uWW$ was defined, assuming 
that $\kk=\CC$; the definition in \cite{\KL} was in terms of the Lie algebra of $G$ with scalars extended to the 
power series field $\CC((\e))$. (The idea that a relationship between $\uuG$ and $\uWW$ might exist appeared in
Carter's paper \cite{\CA}.) In this paper, developing an idea in \cite{\COR}, we define a surjective map
$\Ph:\uWW@>>>\uuG$. Our definition of $\Ph$ is not in terms of the Lie algebra but in terms of the group and it
works in any characteristic (but for the purposes of this introduction we assume that $p$ is not a bad prime for 
$G$). More precisely, we look at the intersection of a Bruhat double coset of $G$ with various unipotent conjugacy
classes and we select the minimal unipotent class which gives a nonempty intersection. (We assume that the Bruhat
double coset corresponds to a Weyl group element which has minimal length in its conjugacy class.) The fact that such a procedure might work is suggested by
the statement in Steinberg \cite{\ST, 8.8} that the Bruhat  double 
coset corresponding to a Coxeter element of minimal length intersects exactly one unipotent class (the regular one), by the result in Kawanaka \cite{\KA} 
that the regular unipotent class of $G$ intersects
every Bruhat double coset, and by the examples in rank $\le3$ given in \cite{\COR}. But the fact that the procedure
actually works is miraculous. In this paper the proof is given separately for classical groups; for exceptional 
groups the desired result is reduced, using the representation theory of reductive groups over a finite field, 
to a computer calculation, see 1.2. I thank Gongqin Li for doing the programming involved in the calculation.

\subhead 0.2\endsubhead
Here is some notation that we use in this paper. Let $Z_G$ be the centre of $G$.
Let $\cb$ be the variety of Borel subgroups of $G$. Let $\WW$ be a set indexing the set of orbits of $G$
acting on $\cb\T\cb$ by $g:(B,B')\m(gBg\i,gB'g\i)$. For $w\in\WW$ we write $\co_w$ for the corresponding $G$-orbit
in $\cb\T\cb$. Define $\ul:\WW@>>>\NN$ by $\ul(w)=\dim\co_w-\dim\cb$. Let $S=\{s\in\WW;\ul(s)=1\}$. There is a 
unique group structure on $\WW$ such that $s^2=1$ for all $s\in S$ and such that
$$\align&w\in\WW,w'\in\WW,(B_1,B_2)\in\co_w,(B_2,B_3)\in\co_{w'},\ul(ww')=\ul(w)+\ul(w')\imp\\&
(B_1,B_3)\in\co_{ww'}.\endalign$$
Then $\WW,S$ is a finite Coxeter group with length function $\ul$ (the Weyl group of $G$). Let $\uWW$ be the set 
of conjugacy classes in $\WW$. For any $C\in\uWW$ let $d_C=\min_{w\in C}\ul(w)$ and let 
$C_{min}=\{w\in C;\ul(w)=d_C\}$. For any $w\in\WW$ let 
$$\fB_w=\{(g,B)\in G\T\cb;(B,gBg\i)\in\co_w\}.$$
(This variety enters in an essential way in the definition of character sheaves on $G$.) We have a partition 
$\fB_w=\sqc_{\g\in\uG}\fB^\g_w$ where
$$\fB^\g_w=\{(g,B)\in\fB_w;g\in\g\}.$$
Note that $G$ and $G_{ad}:=G/Z_G$ act on $\fB_w$ and on $\fB_w^\g$ (for $\g\in\uG$) by 
$x:(g,B)\m(xgx\i,xBx\i)$, $xZ_G:(g,B)\m(xgx\i,xBx\i)$. 
For $\g\in\uG,C\in\uWW$ we write $C\dsv\g$ if $\fB^\g_w\ne\em$ for some/any $w\in C_{min}$. (The equivalence of 
some/any follows from 1.2(a), using \cite{\GP, 8.2.6(b)}.) For $\g\in\uG$ we denote by $\bar\g$ the closure of 
$\g$ in $G$.

For any $J\sub S$ let $\WW_J$ be the subgroup of $\WW$ generated by $J$. We say that $C\in\uWW$ is {\it 
elliptic} if $C\cap\WW_J=\em$ for any $J\sneq S$. Let 
$$\uWW_{el}=\{C\in\uWW;C\text{ elliptic}\}.$$   
If $P$ is a parabolic subgroup of $G$ there is a unique subset $J\sub S$ (said to be the type of $P$) such that 
$$\{w\in\WW;(B,B')\in\co_w\text{ for some }B\sub P,B'\sub P\}=\WW_J.$$
For an integer $\s$ we define $\k_\s\in\{0,1\}$ by $\s=\k_\s\mod2$. For two integers $a,b$ we set 
$[a,b]=\{c\in\ZZ;a\le c\le b\}$. The cardinal of a finite set $X$ is denoted by $|X|$ or by $\sha(X)$.
For $g\in G$, $Z(g)$ denotes the centralizer of $g$ in $G$.
Let $C_{cox}$ be the conjugacy class in $\WW$ that contains the Coxeter elements.
For any parabolic subgroup $P$ of $G$ let $U_P$ be the unipotent radical of $P$.

\subhead 0.3\endsubhead
Let $C\in\uWW$. Consider the following property:

$\Pi_C$. {\it There exists $\g\in\uuG$ such that $C\dsv\g$ and such that if $\g'\in\uuG$ and $C\dsv\g'$ then 
$\g\sub\bar\g'$.}
\nl
Note that if $\Pi_C$ holds then $\g$ is uniquely determined; we denote it by $\g_C$.

We state our main result.

\proclaim{Theorem 0.4} Assume that $p$ is not a bad prime for $G$. Then

(i) $\Pi_C$ holds for any $C\in\uWW$;

(ii) the map $\uWW@>>>\uuG$, $C\m\g_C$ is surjective. 
\endproclaim

\subhead 0.5\endsubhead
Recall that $\g\in\uuG$ is {\it distinguished} if for some/any $g\in G$, $g$ is not contained in a Levi subgroup 
of a proper parabolic subgroup of $G$. In 1.1 it is shown how Theorem 0.4 can be deduced from the following 
result.

\proclaim{Proposition 0.6} Assume that $p$ is not a bad prime for $G$. Then 

(i) $\Pi_C$ holds for any $C\in\uWW_{el}$;

(ii) the map $\uWW_{el}@>>>\uuG$, 
$C\m\g_C$ is injective and its image contains all distinguished unipotent classes of $G$.
\endproclaim

The following result provides an alternative definition for the the map in 0.6(ii).

\proclaim{Theorem 0.7} Assume that $p$ is not a bad prime for $G$. Let $C\in\uWW_{el}$. Let $w\in C_{min}$.

(a) If $\g\in\uG$ and $\fB_w^\g\ne\em$ then $\dim Z(g)/Z_G\le d_C$ for some/all $g\in\g$.

(b) There is a unique unipotent class $\g$ in $G$ such that $\fB_w^\g\ne\em$ and \lb
$\dim Z(g)/Z_G=d_C$ for some/all $g\in\g$.

(c) The class $\g$ in (b) depends only on $C$, not on $w$. It coincides with $\g_C$ in 0.6(ii).
\endproclaim
This follows from results in \S5. 

\subhead 0.8\endsubhead
This paper is organized as follows.
Section 1 contains some preparatory material. In Section 2 we define a particular class of reduced decompositions
for certain elliptic elements of $\WW$. To such a decomposition we attach a unipotent element in $G$. We
study this element in several cases arising from classical groups. This provides one of the ingredients in the
proof of 0.6 for classical groups. (It might also provide an alternative definition for our map $\uWW@>>>\uuG$, 
see the conjecture in 4.7.) In Section 3 we complete the proof of 0.6 for classical groups. In Section 4 we extend
our results to arbitrary characteristic. In 4.3 we give an explicit description of the restriction of the map 
$\Ph:\uWW@>>>\uuG$ to $\uWW_{el}$ for various almost simple $G$. In \S5 we associate to any $C\in\uWW_{el}$ a 
collection of conjugacy classes in $G$, said to be $C$-small classes: the conjugacy classes $\g\in\uG$ of 
minimum dimension such that $C\dsv\g$; we also verify 0.7.

\subhead 0.8\endsubhead
For earlier work on the intersection of Bruhat double cosets with conjugacy classes in $G$ see \cite{\EG}. (I 
thank Jiang-Hua Lu for this reference.)

\head Contents\endhead
1. Preliminaries.

2. Excellent decompositions and unipotent elements.

3. Isometry groups.

4. Basic unipotent classes.

5. $C$-small classes.

\head 1. Preliminaries\endhead
\subhead 1.1\endsubhead
We show how 0.4 can be proved assuming that 0.6 holds when $G$ is replaced by any Levi subgroup of a parabolic
subgroup of $G$. Let $C\in\uWW$. If $C\in\uWW_{el}$ then the result follows from our assumption. We now assume 
that $C$ is not elliptic. We can find $J\sneq S$ and an elliptic conjugacy class $D$ of the Weyl group 
$\WW_J$ such that $D=C\cap\WW_J$. Let $P$ be a parabolic subgroup of $G$ of type $J$. Let $L$ be a Levi subgroup 
of $P$. Let $\g_D$ be the unipotent class of $L$ associated to $D$ by 0.6(i) with $G,\WW$ replaced by $L,\WW_J$. 
Let $\g$ be the unipotent class of $G$ containing $\g_D$. Let $g\in\g$, $w\in D$. Note that some $G$-conjugate 
$g'$ of $g$ is contained in $L$. We can find Borel subgroups $B,B'$ of $P$ such that $(B,B')\in\co_w$, 
$B'=g'Bg'{}\i$. Thus $\fB_w^\g\ne\em$. Now let $\g'$ be a unipotent class of $G$ such that $\fB_{w'}^{\g'}\ne\em$ 
for some $w'\in C_{min}$. We have $C_{min}\cap D\ne\em$ (see \cite{\GP, 3.1.14}) hence we can assume that 
$w'\in D$. We can find $(B,B')\in\co_{w'}$ and $g'\in\g'$ such that $B'=g'Bg'{}\i$. Replacing $B,B',g'$ by 
$xBx\i,xB'x\i,xg'x\i$ for some $x\in G$ we see that we can assume that $B\sub P$ and then we automatically have 
$B'\sub P$ that is $g'Bg'{}\i\sub P$. We have also $g'Bg'{}\i\sub g'Pg'{}\i$ hence $g'Pg'{}\i=P$ that is, 
$g'\in P$. We have $g'=g'_1v$ where $g'_1\in L$ is unipotent and $v\in U_P$. 
We can find a one parameter subgroup $\l:\kk^*@>>>Z_L$ such that $\l(t)v\l(t\i)$ converges 
to $1$ when $t\in\kk^*$ converges to $0$. Then $\l(t)g'\l(t)\i=g'_1\l(t)v\l(t)\i$ converges to $g'_1$ when 
$t\in\kk^*$ converges to $0$. Thus $g'_1$ is contained in the closure of $\g'$. Hence the $L$-conjugacy class
of $g'_1$ is contained in the closure of $\g'$. Note also that $B'=g'_1Bg'_1{}\i$. Using the definition of $\g_D$
we see that $\g_D$ is contained in the closure of the $L$-conjugacy class of $g'_1$. Hence $\g_D$ is contained in
the closure of $\g'$ and $\g$ is contained in the closure of $\g'$. We see that $\g$ has the property stated in 
$\Pi_C$. This proves 0.4(i) (assuming 0.6(i)).

The previous argument shows that $\g_C$ is the unipotent class of $G$ containing the unipotent class $\g_D$ of
$L$. Thus $C\m\g_C$ is determined in a simple way from the knowledge of the maps $D\m\g_D$ in 0.6 corresponding to
various $L$ as above.

Now let $\g\in\uuG$. We can find a parabolic subgroup $P$ of $G$ with Levi subgroup $L$ and
a distinguished unipotent class $\g_1$ of $L$ such that $\g_1\sub\g$. Let $J$ be the subset $S$ such that $P$ is 
of type $J$. By 0.6(ii) we can find an elliptic conjugacy class $D$ of $\WW_J$ such that $\g_D=\g_1$ (where 
$\g_D$ is defined in terms of $L,D$). Let $C$ be the conjugacy class in $\WW$ that contains $D$. By the arguments
above we have $\g_C=\g$. This proves 0.4(ii) (assuming 0.6).

\subhead 1.2\endsubhead
To prove Proposition 0.6 we can assume that $G$ is almost simple. Moreover for each isogeny class of almost simple
groups it is enough to prove 0.6 for one group in the isogeny class and 0.6 will be automatically true for the 
other groups in the isogeny class.

Note that if $C=C_{cox}$ (recall that $C\in\uWW_{el}$) then 
$\Pi_C$ follows from a statement in \cite{\ST, 8.8}; in this case $\g_C$ is the regular unipotent class. If $G$ 
is almost simple of type $A_n$ then $C$ as above is the only element of $\uWW_{el}$ and the only distinguished 
unipotent class is the regular one so that in this case 0.6 follows.

If $G$ is almost simple of type $B_n,C_n$ or $D_n$ then we can assume that $G$ is as in 1.3. The proof of 0.6 in 
these cases is given in 3.7-3.9. 

In the remainder of this subsection we assume that $\kk$ is an algebraic closure of a finite field $\FF_q$ with 
$q$ elements. We choose an $\FF_q$-split rational structure on $G$ with Frobenius map $F:G@>>>G$. Now $F$ induces
a morphism $\cb@>>>\cb$ denoted again by $F$. Note that the finite group $G^F$ acts transitively on the finite 
set $\cb^F$ (the upper script denotes the set of fixed points). 
Hence $G^F$ acts naturally on the $\bbq$-vector space $\cf$ of functions $\cb^F@>>>\bbq$. (Here $l$ is a fixed
prime number such that $l\ne0$ in $\kk$.) For any $w\in\WW$ we denote by $T_w:\cf@>>>\cf$ the linear map $f\m f'$
where $f'(B)=\sum_{B'\in\cb^F;(B,B')\in\co_w}f(B')$. Let $\ch_q$ be the subspace of $\End(\cf)$ spanned by
$T_w(w\in\WW$); this is a subalgebra of $\End(\cf)$ and the irreducible $\ch_q$-modules (up to
isomorphism) are in natural bijection $E_q\lra E$ with $\Irr\WW$, the set of irreducible $\WW$-modules over $\bbq$
(up to isomorphism) once $\sqrt{q}$ has been chosen. Moreover we have a canonical decomposition 
$\cf=\op_{E\in\Irr\WW}E_q\ot\r_E$ (as a $(\ch_q,G^F)$-module) where $\r_E$ is an irreducible representation of 
$G^F$. Now let $\g$ be an $F$-stable $G$-conjugacy class in $G$. Then $\fB^\g_w$ has a natural Frobenius map
$(g,B)\m(F(g),F(B))$ denoted again by $F$. We compute the number of fixed points of $F:\fB^\g_w@>>>\fB^\g_w$:
$$\align&|(\fB^\g_w)^F|=\sum_{g\in\g^F}\sha(B\in\cb^F;(B,gBg\i)\in\co_w\}\\&=
\sum_{g\in\g^F}\tr(gT_w:\cf@>>>\cf)=\sum_{g\in\g^F}\sum_{E\in\Irr\WW}\tr(T_w,E_q)\tr(g,\r_E).\tag a\endalign$$
For any $y\in\WW$ let $R^\th(y)$ be the virtual representation of $G^F$ defined in \cite{\DL, 1.9} ($\th$ as in
\cite{\DL, 1.8}). We have
$\r_E=|\WW|\i\sum_{y\in\WW}(\r_E:R^1(y))R^1(y)+\x_E$ where $\x_E$ is a virtual 
representation of $G^F$ orthogonal to each $R^\th(y)$ and 
$(\r_E:R^1(y))$ denotes multiplicity. Using the equality
$$\sum_{g\in\g^F}\tr(g,\x_E)=0\tag b$$
(verified below) we deduce that
$$|(\fB^\g_w)^F|=\sum_{g\in\g^F}\sum_{E\in\Irr\WW}\tr(T_w,E_q)|\WW|\i\sum_{y\in\WW}(\r_E:R^1(y))\tr(g,R^1(y)).$$
For any $E'\in\Irr\WW$ we set $R_{E'}=|\WW|\i\sum_{y'\in\WW}\tr(y',E')R^1(y')$ so that
$R^1(y)=\sum_{E'\in\Irr\WW}\tr(y,E')R_{E'}$. We have
$$\align&|(\fB^\g_w)^F|\\&=\sum_{g\in\g^F}\sum_{E,E',E''\in\Irr\WW}
\tr(T_w,E_q)|\WW|\i\sum_{y\in\WW}\tr(y,E')\tr(y,E'')(\r_E:R_{E'})\tr(g,R_{E''}).\endalign$$
Hence 
$$|(\fB^\g_w)^F|=|\WW|\i\sum_{E,E'\in\Irr\WW,y'\in\WW,g\in\g^F}\tr(T_w,E_q)(\r_E:R_{E'})\tr(y',E')\tr(g,R^1(y')).
\tag c$$
We now verify (b). Let $A$ be the vector space of $G^F$-invariant invariant functions $G^F@>>>\bbq$. Let $A_0$ be
the subspace of $A$ spanned by the functions $f^\th(y)$ (the character of $R^\th(y)$) for various $y,\th$. Let 
$\c\in A$ be the characteristic function of the subset $\g^F$ of $G^F$. We must show that if $f\in A$ is 
orthogonal to $A_0$ then it is orthogonal to $\c$. It is enough to show that $\c\in A_0$. If $G$ is a classical 
group or if $p$ is not a bad prime for $G$, this follows from results in \cite{\CS,\GF}; in the general case it is
proved by M. Geck \cite{\GE}, using results in \cite{\CS,\GF,\SH}.

Now assume that $G$ is adjoint of exceptional type, that $p$ is not a bad prime for $G$, that $\g\in\uuG$ and 
that $w\in C_{min}$ where $C\in\uWW_{el}$. We also assume that $q-1$ is sufficiently divisible. Then (c) becomes
$$|(\fB^\g_w)^F|=|\WW|\i\sum A_{E,C}\ph_{E,E'}a_{E',C'}d_{C',C''}Q(C'',\g_0)D(\g_0,\g'_0)P(\g'_0,\g)$$
where the sum is taken over all $E,E'$ in $\Irr\WW$, $C',C''$ in $\uWW$, $\g_0,\g'_0$ in $\uuG_0$ and the
notation is as follows. 

$\uuG_0$ is the set of $G^F$-conjugacy classes of unipotent elements in $G^F$. For $\g_0\in\uuG_0,\g'\in\uuG$ we 
set $D_{\g_0,\g'}=|\g_0|$ if $\g_0\sub\g'$ and $D_{\g_0,\g'}=0$ if $\g_0\not\sub\g'$. For $C',C''\in\uWW$ we set 
$d_{C',C''}=|C'|$ if $C'=C''$ and $d_{C',C''}=0$ if $C'\ne C''$. For $C'\in\uWW,E\in\Irr\WW$ we set 
$A_{E,C'}=\tr(T_z,E_q)$, $a_{E,C'}=\tr(z,E)$ where $z\in C'_{min}$. (Note that $A_{E,C'}$ is well defined by 
\cite{\GP, 8.2.6(b)}.) For $E,E'\in\Irr\WW$ let $\ph_{E,E'}=(\r_E:R_{E'})$. For $\g_0\in\uuG_0$ and $C\in\uWW$ 
let $Q_{C,\g_0}=\tr(g,R^1(y))$ where $g\in\g_0,y\in C$.

Thus $|(\fB^\g_w)^F|$ is $|\WW|\i$ times the $(C,\g)$ entry of the matrix which is the product of matrices
$${}^t(A_{E,C})(\ph_{E,E'})(a_{E',C'})(d_{C',C''})(Q_{C'',\g_0})(D_{\g_0,\g}).$$
Each of these matrices is explicitly known. The matrix $(A_{E,C'})$ is known from the works of Geck and 
Geck-Michel (see \cite{\GP, 11.5.11}) and is available through the CHEVIE package; the matrix $(d_{C',C''})$ is 
available from the same source. The matrix $(a_{E,C'})$ is the specialization $q=1$ of $(A_{E,C'})$. The matrix 
$\ph_{E,E'}$ has as entries the coefficients of the "nonabelian Fourier transform" in \cite{\OR, 4.15}. The matrix
$(Q_{C'',\g_0})$ is the matrix of Green functions, known from the work of Shoji and Beynon-Spaltenstein. I thank 
Frank L\"ubeck for providing tables of Green functions in GAP-format and instructions on how to use them; these 
tables can now be found at \cite{\LUE}; the matrix $(D_{\g_0,\g'})$ can be extracted from the same source. Thus 
$|(\fB^\g_w)^F|$ can be obtained by calculating the product of six (large) explicitly known matrices. The 
calculation was done using the CHEVIE package, see \cite{\GH}. It turns out that $|(\fB^\g_w)^F|$ is a polynomial
in $q$ with integer coefficients. Note that $\fB^\g_w\ne\em$ if and only if $|(\fB^\g_w)^F|\ne0$ for sufficiently
large $q$. Thus the condition that $C\dsv\g$ can be tested. This can be used to check that 0.6 holds in our case.
(This method is a simplification of the method in \cite{\COR, 1.5}.) 

From the explicit calculations above we see that the following hold when $C$ is elliptic:

(d) {\it If $\g=\g_C,w\in C_{min}$, then $|(\fB^\g_w)^F|/|G^F|$ is a polynomial in $q$ with constant term $1$. If
$C\dsv\g,w\in C_{min}$ but $\g\ne\g_C$, then $|(\fB^\g_w)^F|/|G^F|$ is a polynomial in $q$ with costant term $0$.
If $w\in C_{min}$, the sum $\sum_{\g\in\uuG}|(\fB^\g_w)^F|/|G^F|$ is a palindromic polynomial in $q$ of the form
$1+\do+q^{\ul(w)-r}$ ($r$ is the rank of $G_{ad}$); the constant term $1$ comes from $\g=\g_C$ and the highest
term $q^{\ul(w)-r}$ comes from the regular unipotent class.}
\nl
We now see that 0.6 holds. The correspondence $C\m\g_C$ for $C\in\uWW_{el}$ is described explicitly in 4.3.

We expect that (d) also holds for classical types.

\subhead 1.3\endsubhead
Let $V$ be a $\kk$-vector space of finite dimension $\nn\ge3$. We set $\k=\k_\nn$. Let $n=(\nn-\k)/2$. Assume that
$V$ has a fixed bilinear form $(,):V\T V@>>>\kk$ and a fixed quadratic form $Q:V@>>>\kk$ such that (i) or (ii) 
below holds:

(i) $Q=0$, $(x,x)=0$ for all $x\in V$, $V^\pe=0$;

(ii) $Q\ne0$, $(x,y)=Q(x+y)-Q(x)-Q(y)$ for $x,y\in V$, $Q:V^\pe@>>>\kk$ is injective.
\nl
Here, for any subspace $V'$ of $V$ we set $V'{}^\pe=\{x\in V;(x,V')=0\}$. In case (ii) it follows that $V^\pe=0$ 
unless $\k=1$ and $p=2$ in which case $\dim V^\pe=1$. An element 
$g\in GL(V)$ is said to be an isometry if $(gx,gy)=(x,y)$ for all $x,y\in V$ and $Q(gx)=Q(x)$ for all $x\in V$.
Let $Is(V)$ be the group of all isometries of $V$ (a closed subgroup of $GL(V)$). A subspace $V'$ of $V$ is said 
to be isotropic if $(,)$ and $Q$ are zero on $V'$. Let $\cf$ be the set of all sequences 
$V_*=(0=V_0\sub V_1\sub V_2\sub\do\sub V_\nn=V)$ of subspaces of $V$ such that $\dim V_i=i$ for $i\in[0,\nn]$, 
$Q|_{V_i}=0$ and $V_i^\pe=V_{\nn-i}$ for all $i\in[0,n]$. (For such $V_*$, $V_i$ is an isotropic subspace for 
$i\in[0,n]$). Now $Is(V)$ acts naturally (transitively) on $\cf$. 

In the remainder of this section we assume that $G$ is the identity component of $Is(V)$. 

\subhead 1.4\endsubhead
Let $W$ be the group of permutations of $[1,\nn]$ which commute with the involution $i\m\nn-i+1$ of $[1,\nn]$. (In
particular, if $\k=1$ then any permutation in $W$ fixes $n+1$.) Let $V_*,V'_*$ be two sequences in $\cf$. As in 
\cite{\HEC, 0.4} we define a permutation $a_{V_*,V'_*}:i\m a_i$ of $[1,\nn]$ as follows. For 
$i\in[0,\nn],j\in[1,\nn]$ we set $d_{ij}=\dim(V'_i\cap V_j)/(V'_i\cap V_{j-1})\in\{0,1\}$. For $i\in[0,\nn]$ we 
set $X_i=\{j\in[1,\nn];d_{ij}=1\}$. We have $\em=X_0\sub X_1\sub X_2\sub\do\sub X_{\nn}=[1,\nn]$ and for 
$i\in[1,\nn]$ there is a unique $a_i\in[1,\nn]$ such that $X_i=X_{i-1}\sqc\{a_i\}$. Then $i\m a_i$ is the 
required permutation of $[1,\nn]$. It belongs to $W$. Moreover 

(a) {\it $(V_*,V'_*)\m a_{V_*,V'_*}$ defines a bijection from the set of $Is(V)$-orbits on $\cf\T\cf$ (for the 
diagonal action) to $W$.}
\nl
When $\k=0,Q\ne0$ we define $W'$ as the group of even permutations in $W$ (a subgroup of index $2$ of $W$). For 
$i\in[1,n-1]$ define $s_i\in W$ as a product of two transpositions $i\lra i+1$, $\nn+1-i\lra\nn-i$ (all other
entries go to themselves); define $s_n\in W$ to be the transposition $n\lra\nn-n+1$ (all other entries go to 
themselves). Then $(W,\{s_i;i\in[1,n]\})$ is a Weyl group of type $B_n$. If $\k=0,Q\ne0$, we have $s_i\in W'$ for
$i\in[1,n-1]$ and we set $\ts_i=s_ns_is_n\in W'$ for $i\in[1,n-1]$; we have $\ts_i=s_i$ if $i<n-1$ and 
$(W',\{s_1,s_2,\do,s_{n-1},\ts_{n-1}\})$ is a Weyl group of type $D_n$.

\subhead 1.5\endsubhead 
For any $V_*\in\cf$ we set $B_{V_*}=\{g\in G;gV_*=V_*\}$, a Borel subgroup of $G$. If $(1-\k)Q=0$ then 
$V_*\m B_{V_*}$ is an isomorphism $\cf@>\si>>\cb$; for $w\in W$ and 
$V_*,V'_*$ in $\cf$ such that $a_{V_*,V'_*}=w$ we have $(B_{V_*},B_{V'_*})\in\co_{f(w)}$ for a well defined 
element $f(w)\in\WW$ and $w\m f(w)$ is an isomorphism $W@>\si>>\WW$ (as Coxeter groups) by which these two groups
are identified. If $(1-\k)Q\ne0$, $G$ has 
two orbits on $\cf$. Let $\cf'$ be one of these orbits. Then $V_*\m B_{V_*}$ is an isomorphism $\cf'@>\si>>\cb$; 
for $w\in W'$ and $V_*,V'_*$ in $\cf$ such that $a_{V_*,V'_*}=w$ and $V_*\in\cf'$ we have $V'_*\in\cf'$ and 
$(B_{V_*},B_{V'_*})\in\co_{f(w)}$ for a well defined element $f(w)\in\WW$. Moreover $w\m f(w)$ is an isomorphism 
$W'@>\si>>\WW$ (as Coxeter groups) by which these two groups are identified.

\subhead 1.6\endsubhead
Let $\cp_n$ be the set of sequences $p_1\ge p_2\ge\do\ge p_\s$ of integers $\ge1$ such that $p_1+p_2+\do+p_\s=n$.
Let $\cp_n^+$ be the set of sequences $p_1\ge p_2\ge\do\ge p_\s$ in $\cp_n$ such that $\k_\s=0$. For any 
$r\in[1,\s]$ we set $p_{\le r}=\sum_{r'\in[1,r]}p_{r'}$, $p_{<r}=\sum_{r'\in[1,r-1]}p_{r'}$,
$p_{>r}=\sum_{r'\in[r+1,\s]}p_{r'}$. For $p_*\in\cp_n$ let $w_{p_*}\in W$ be the permutation of $[1,\nn]$ given by
$$\align&1\m2\m\do\m p_1\m\nn\m\nn-1\m\do\m\nn-p_1\m1,\\&
p_1+1\m p_1+2\m\do\m p_1+p_2\m\nn-p_1-1\m\nn-p_1-2\m\do \\&  \m\nn-p_1-p_2\m p_1+1,\\&\do\\&
p_{<\s}+1\m p_{<\s}+2\m\do\m p_{<\s}+p_\s\m\nn-p_{<s}-1\\&\m\nn-p_{<\s}-2\m\do\m\nn-p_{<s}-p_\s\m p_{<\s}+1,\\&
\text{ and, if }\k=1,\\&p_{n+1}\m p_{n+1}.\tag a\endalign$$
Let $C_{p_*}$ be the conjugacy class of $w_{p_*}$ in $W$. If $(1-\k)Q=0$, then $w_{p_*}\in\uWW_{el}$ and 
$p_*\m C_{p_*}$ is a bijection $\cp_n@>\si>>\uWW_{el}$. If $(1-\k)Q\ne0$ and $p_*\in\cp_n^+$ then 
$w_{p_*}\in\uWW_{el}$; we denote by $C'_{p_*}$ the conjugacy class of $w_{p_*}$ in $W'$. Then $p_*\m C'_{p_*}$ is
a bijection $\cp_n^+@>\si>>\uWW_{el}$. 

For any $p_1\ge p_2\ge\do\ge p_\s$ in $\cp_n$ we define a function $\ps:[1,\s]@>>>\{-1,0,1\}$ as follows.

(i) If $t\in[1,\s]$ is odd and $p_t<p_x$ for any $x\in[1,t-1]$ then $\ps(t)=1$;

(ii) if $t\in[1,\s]$ is even and $p_x<p_t$ for any $x\in[t+1,\s]$, then $\ps(t)=-1$;

(iii) for all other $t\in[1,\s]$ we have $\ps(t)=0$.
\nl
For any integer $a$ such that $1\le 2a<2a+1\le\s$ we have $\ps(2a)+\ps(2a+1)=0$. (Indeed if $p_{2a}>p_{2a+1}$ 
then $\ps(2a)+\ps(2a+1)=-1+1=0$; if $p_{2a}=p_{2a+1}$ then $\ps(2a)+\ps(2a+1)=0+0=0$.) Also $\ps(1)=1$. Hence 

{\it if $h\in[1,\s]$ is odd, then $\sum_{r\in[1,h]}\ps(r)=1$;}

{\it if $h\in[1,\s]$ is even, then $\sum_{r\in[1,h]}\ps(r)=1+\ps(h)$.}
\nl
We have $\ps(\s)=-1$ if $\s$ is even. Hence
$$\sum_{t\in[1,\s]}\ps(t)=\k_\s.\tag b$$

\head 2. Excellent decompositions and unipotent elements \endhead
\subhead 2.1\endsubhead
Let $C\in\uWW_{el}$ and let $w\in C$. Let $\r=|S|$. An {\it excellent decomposition} of $w$ is a sequence
$$\align&s_1^1,s_2^1,\do,s_{q_1}^1,s_{q_1+1}^1,s_{q_1}^1,\do,s_2^1,s_1^1,
s_1^2,s_2^2\do,s_{q_2}^2,s_{q_2+1}^2,s_{q_2}^2,\do,s_2^2,s_1^2,\do,\\&
s_1^\r,s_2^\r,\do,s_{q_\r}^\r,s_{q_\r+1}^\r,s_{q_\r}^\r,\do,s_2^\r,s_1^\r\endalign$$
in $S$ (the upper scripts are not powers) consisting of $\sum_{k\in[1,\r]}(2q_k+1)=\ul(w)$ terms, such that
$$\align&w=(s_1^1s_2^1\do s_{q_1}^1s_{q_1+1}^1s_{q_1}^1\do s_2^1s_1^1)
(s_1^2s_2^2\do s_{q_2}^2s_{q_2+1}^2s_{q_2}^2\do s_2^2s_1^2)\do\\&
(s_1^\r s_2^\r\do s_{q_\r}^\r s_{q_\r+1}^\r s_{q_\r}^\r\do s_2^\r s_1^\r).\tag a\endalign$$
We sometime refer to (a) as an excellent decomposition of $w$. It is a reduced expression of a special kind for 
$w$. It appears that

(b) {\it for any $C\in\uWW_{el}$, at least one element $w\in C_{min}$ admits an excellent decomposition.}
\nl
For example if $C=C_{cox}$ then for any $w\in C_{min}$ any reduced expression
of $w$ is an excellent decomposition. In particular (b) holds when $G$ is almost simple of type $A_n$. When $G$ is
simple of type $G_2$, the excellent decompositions $(1)(2),(121)(2),(12121)(2)$ account for the $3$ elliptic
conjugacy classes in $\WW$ with $S=\{1,2\}$. In type $F_4$, the excellent decompositions 
$$\align&(1)(2)(3)(4);\qua (1)(232)(3)(4);\qua (121)(323)(4)(3);\qua (1)(2)(3234323)(4);\\&(4)(3)(2321232)(1);\qua
(12321)(23432)(3)(4);\qua(2)(12321)(3234323)(4);\\&(2324312134232)(3)(1)(4);\qua(432134232431234)(12321)(232)(3)
\endalign$$
(notation of \cite{\GP, p.407}) account for the $9$ elliptic conjugacy classes in $\WW$; in type $E_6$, the 
excellent decompositions 
$$\align&(1)(2)(3)(4)(5)(6);\qua (1)(3)(4)(2)(454)(6);\qua (1)(3)(4)(2345432)(6)(5);\\&(1)(2)(3)(432454234)(5)(6);
\qua(4354132456542314534)(2)(1)(3)(5)(6)\endalign$$
(notation of \cite{\GP, p.407}) account for the $5$ elliptic conjugacy classes in $\WW$. Note that some 
elements in $C_{min}$ might not admit an excellent decomposition (example: the element $324312$ in type $F_4$). 
In 2.2, 2.3 we will verify (b) for $G$ of type $B_n,C_n,D_n$.

\subhead 2.2\endsubhead
Assume that notation is as in 1.3 and that $(1-\k)Q=0$. Let $p_*=(p_1\ge p_2\ge\do\ge p_\s)\in\cp_n$. The
following is an excellent decomposition of $w=w_{p_*}\i$ (see 1.6(a)) in $\WW$:
$$\align&w=(s_n)(s_{n-1})\do(s_{n-p_\s+1})\T\\&
(s_{n-p_\s}\do s_{n-1}s_ns_{n-1}\do s_{n-p_\s})(s_{n-p_\s-1})(s_{n-p_\s-2})\do(s_{n-p_\s-p_{s-1}+1})\T\\&
(s_{n-p_\s-p_{\s-1}}\do s_{n-1}s_ns_{n-1}\do s_{n-p_\s-p_{\s-1}})(s_{n-p_\s-p_{\s-1}-1})(s_{n-p_\s-p_{\s-1}-2})\\&
\do(s_{n-p_\s-p_{\s-1}-p_{\s-2}})\T\\&
\do\\&(s_{n-p_\s-\do-p_2}\do s_{n-1}s_ns_{n-1}\do s_{n-p_\s-\do-p_2})(s_{n-p_\s-\do-p_2-1})\\&
(s_{n-p_\s-\do-p_2-2})\do(s_{n-p_\s-\do-p_1+1}).\tag a\endalign$$
Note that $\ul(w)=2\sum_{v=1}^{\s-1}vp_{v+1}+n$. Thus $w$ has minimal length in its conjugacy class in $\WW$ (see
\cite{\GP, 3.4}). We see that 2.1(b) holds for $G$ of type $B_n$ or $C_n$.

In the case where $p_*\in\cp^+_n$ we define another excellent decomposition of $w=w_{p_*}\i$ (only the parantheses
differ from the previous excellent decomposition):
$$\align&w=(s_ns_{n-1}\do s_{n-p_\s}\do s_{n-1}s_n)(s_{n-1})(s_{n-2})\do(s_{n-p_\s-p_{\s-1}+1})\T\\&
(s_{n-p_\s-p_{\s-1}}\do s_{n-1}s_ns_{n-1}\do s_{n-p_\s-p_{\s-1}-p_{\s-2}}\do s_{n-1}s_ns_{n-1}\do\\& 
s_{n-p_\s-p_{\s-1}})(s_{n-p_\s-p_{\s-1}-1})(s_{n-p_\s-p_{\s-1}-2})\do(s_{n-p_\s-p_{\s-1}-p_{\s-2}-p_{\s-3}+1})
\T\\&\do\\&(s_{n-p_\s-\do-p_3}\do s_{n-1}s_ns_{n-1}\do s_{n-p_\s-\do-p_2}\do s_{n-1}s_ns_{n-1}\do \\&
s_{n-p_\s-\do-p_3})(s_{n-p_\s-\do-p_3-1})(s_{n-p_\s-\do-p_3-2})\do(s_{n-p_\s-\do-p_1+1}).\tag b\endalign$$

\subhead 2.3\endsubhead
Assume that notation is as in 1.3 and that $(1-\k)Q\ne0$. Let $p_*=(p_1\ge p_2\ge\do\ge p_\s)\in\cp^+_n$. The 
excellent decomposition 2.2(b) in $W$ gives rise to an excellent decomposition of $w=w_{p_*}\i$ in $W'=\WW$:
$$\align&w=(\ts_{n-1}\do\ts_{n-p_\s}\do\ts_{n-1})(s_{n-1})(s_{n-2})\do(s_{n-p_\s-p_{\s-1}+1})\T\\&
(s_{n-p_\s-p_{\s-1}}\do s_{n-1}\ts_{n-1}\do\ts_{n-p_\s-p_{\s-1}-p_{\s-2}}\do\ts_{n-1}s_{n-1}\do \\&
s_{n-p_\s-p_{\s-1}})(s_{n-p_\s-p_{\s-1}-1})(s_{n-p_\s-p_{\s-1}-2})\do(s_{n-p_\s-p_{\s-1}-p_{\s-2}-p_{\s-3}+1})\T
\\&\do\\&(s_{n-p_\s-\do-p_3}\do s_{n-1}\ts_{n-1}\do\ts_{n-p_\s-\do-p_2}\do\ts_{n-1}s_{n-1}\do \\&
s_{n-p_\s-\do-p_3})(s_{n-p_\s-\do-p_3-1})(s_{n-p_\s-\do-p_3-2})\do(s_{n-p_\s-\do-p_1+1}).\endalign$$
Now the length of $w$ in $W'$ is equal to the length of $w$ in $W$ minus $\s$; hence it is 
$2\sum_{v=1}^{\s-1}vp_{v+1}+n-\s$. Thus $w$ has minimal length in its conjugacy class in $\WW=W'$ (see 
\cite{\GP, 3.4}). We see that 2.1(b) holds for $G$ of type $D_n$.

\subhead 2.4\endsubhead
We return to the general case. We choose a Borel subgroup $B$ of $G$ and a Borel subgroup $B'$ opposed to $B$. Let
$T=B\cap B'$, a maximal torus $T$ of $G$. Let $N(T)$ be the normalizer of $T$ in $G$. Let $U'=U_{B'}$. For any 
$s\in S$ let $P_s$ be the parabolic subgroup of type $\{s\}$ that contains $B$. Note that
$U'_s:=U'\cap P_s$ is isomorphic to $\kk$. Let $t\m y_s(t)$ be an isomorphism of algebraic groups
$\kk@>\si>>U'_s$. Let $\ds\in NT$ be an element such that $(B,\ds B\ds\i)\in\co_s$. 

Let $C\in\uWW_{el}$ and let 2.1(a) be an excellent decomposition of an element $w\in C_{min}$. Let
$c_1,c_2,\do,c_\r$ be elements of $\kk^*$. We set 
$$\align&u_w=(\ds_1^1\ds_2^1\do\ds_{q_1}^1y_{s_{q_1+1}^1}(c_1)(\ds_{q_1}^1)\i\do(\ds_2^1)\i(\ds_1^1)\i)\\&
(\ds_1^2\ds_2^2\do\ds_{q_2}^2y_{s_{q_2+1}^2}(c_2)(\ds_{q_2}^2)\i\do(\ds_2^2)\i(\ds_1^2)\i)\do\\&
(\ds_1^\r\ds_2^\r\do\ds_{q_\r}^\r y_{s_{q_\r+1}^\r}(c_\r)(\ds_{q_\r}^\r)\i\do(\ds_2^\r)\i(\ds_1^\r)\i).\tag a
\endalign$$
We have $u_w\in U'$; thus, $u_w$ is unipotent. Since $y_{s_{q_h+1}^h}(c_h)\in B\ds_{q_h+1}^hB$ and since 2.1(a) is
a reduced expression for $w$, we see (using properties of the Bruhat decomposition) that 
$$(B,u_wBu_w\i)\in\co_w.\tag b$$
Much of the remainder of this section is devoted to computing $u_w$ in some cases arising from classical groups.

\subhead 2.5\endsubhead
Assume that notation is as in 1.3. In the remainder of this section we assume that $(\k-1)Q=0$. We can find (and 
fix) a basis of $V$ consisting of vectors $e_i,e'_i (i\in[1,n])$ and (if $\k=1$) of $e_0$, such that 
$(e_i,e'_j)=\d_{i,j}$ for $i,j\in[1,n]$, $Q(e_i)=Q(e'_i)=0$ for $i\in[1,n]$ and (if $\k=1$), 
$(e_i,e_0)=(e'_i,e_0)$ for $i\in[1,n]$, $(e_0,e_0)=2$, $Q(e_0)=1$.

For $i\in[1,n]$ let $\fV_i$ be the span of $e_1,e_2,\do,e_i$ and let $\fV'_i$ be the span of $e'_1,e'_2,\do,e'_i$.
Let $B$ (resp. $B'$) be the set of all $g\in G$ such that for any $i\in[1,n]$, we have $g\fV_i=\fV_i$ (resp. 
$g\fV'_i=\fV'_i$). Note that $B,B'$ are opposed Borel subgroups of $G$.

We will sometime specify a linear map $V@>>>V$ by indicating its effect on a part of a basis of $V$ with the
understanding that the remaining basis elements are sent to themselves. Thus for $h\in[1,n-1]$ we define 
$y_{s_h}(a)\in GL(V)$ by $[e_h\m e_h+ae_{h+1},e'_{h+1}\m e'_{h+1}-ae'_h]$ ($a\in\kk$). In the case where 
$\k=0,Q=0$ we define $y_{s_n}(a)\in GL(V)$ by $[e_n\m e_n-ae'_n]$, ($a\in\kk$). In the case where $\k=1,Q\ne0$ we 
define $y_{s_n}(a)\in GL(V)$ by $[e_n\m e_n+ae_0-a^2e'_n,e_0\m e_0-2ae'_n]$, $(a\in\kk$). In both cases we have 
$y_{s_h}(a)\in G$. Note that for $s=s_h\in S$, $y_{s_h}:\kk@>>>U'_s$ is as in 2.4.

Let $p_*=(p_1\ge p_2\ge\do\ge p_\s)\in\cp_n$ and let $w=w_{p_*}\i$ with $w_{p_*}$ as in 1.6. We define $u_w$ as in
2.4(a) in terms of $c_1,c_2,\do,c_\r$ in $\kk^*$ and the excellent decomposition 2.2(a) of $w$. Let 
$N=u_w\i-1\in\End(V)$. 

\subhead 2.6\endsubhead
Assume now that $\k=0,Q=0$. From the definitions, for a suitable choice of $c_1,c_2,\do,c_\r$, $u_w\i$ is the 
product over $h\in[1,\s]$ of the linear maps
$$\align&[e_{n-p_{>h}}\m e_{n-p_{>h}}-e'_{n-p_{>h}}+e'_{n-p_{>h}-1}-e'_{n-p_{>h}-2}+\do+(-1)^{p_h}
e'_{n-p_{>h}-p_h+1},
\\&e_{n-p_{>h}-1}\m e_{n-p_{>h}-1}+e_{n-p_{>h}},e_{n-p_{>h}-2}\m e_{n-p_{>h}-2}+e_{n-p_{>h}-1},\\&
\do,e_{n-p_{>h}-p_h+1}\m e_{n-p_{>h}-p_h+1}+e_{n-p_{>h}-p_h+2},\\&
e'_{n-p_{>h}}\m e'_{n-p_{>h}}-e'_{n-p_{>h}-1}+\do+(-1)^{p_h-1}e'_{n-p_{>h}-p_h+1},\\&
e'_{n-p_{>h}-1}\m e'_{n-p_{>h}-1}-e'_{n-p_{>h}-2}+\do+(-1)^{p_h-2}e'_{n-p_{>h}-p_h+1},\\&
\do,e'_{n-p_{>h}-p_h+2}\m e'_{n-p_{>h}-p_h+2}-e'_{n-p_{>h}-p_h+1}].\endalign$$
Hence $u_w\i$ is given by
$$\align&e_i\m e_i+e_{i+1}\text{ if }i\in[1,n],i\n\{p_1,p_1+p_2,\do,p_1+p_2+\do+p_\s\},\\&
e_{p_1+p_2+\do+p_r}\m e_{p_1+p_2+\do+p_r}+\sum_{v\in[1,p_r]}(-1)^ve'_{p_1+p_2+\do+p_r-v+1}
\text{ if }r\in[1,\s],\\&e'_{p_1+p_2+\do+p_r-j}\m\sum_{v\in[1,p_r-j]}(-1)^ve'_{p_1+p_2+\do+p_r-j-v+1}
\text{ if }r\in[1,\s], j\in[0,p_r-2],\\&       
e'_{i_1+i_2+\do+p_{r-1}+1}\m e'_{p_1+p_2+\do+p_{r-1}+1} \text{ if }r\in[1,\s].\endalign$$
We set
$$\align& e^t_j=e_{p_1+\do+p_{t-1}+j}, (t\in[1,\s],j\in[1,p_t]),\\&
e'{}^t_1=\sum_{v\in[1,p_t]}(-1)^ve'_{p_1+p_2+\do+p_t-v+1}, (t\in[1,\s]),\\&
e'{}^t_j=N^{j-1}e'{}^t_1, (t\in[1,\s],j\ge2).\endalign$$
Note that $e'{}^t_j=0$ if $j>p_t$. Clearly for any $t\in[1,\s]$, the elements $e'{}^t_j (j\in[1,p_t])$ span the 
same subspace as the elements $e'_{p_1+p_2+\do+p_t-v+1} (v\in[1,p_t])$. It follows that 
$e^t_j,e'{}^t_j(t\in[1,\s],j\in[1,p_t])$ form a basis of $V$. In this basis the action of $N$ is given by
$$e^t_1\m e^t_2\m\do\m e^t_{p_t}\m e'{}^t_1\m e'{}^t_2\m\do\m e'{}^t_{p_t}\m0 \text{ if }t\in[1,\s].$$
Thus the Jordan blocks of $N:V@>>>V$ have sizes $2p_1,2p_2,\do,2p_\s$.

\subhead 2.7\endsubhead
Assume now that $\k=1,Q\ne0$. From the definitions, for a suitable choice of $c_1,c_2,\do,c_\r$, $u_w\i$ is the 
product over $h\in[1,\s]$ of the linear maps
$$\align&[e_{n-p_{>h}}\m e_{n-p_{>h}}+e_0-e'_{n-p_{>h}}+e'_{n-p_{>h}-1}-e'_{n-p_{>h}-2}+\do+\\&
(-1)^{p_h}e'_{n-p_{>h}-p_h+1},
\\&e_{n-p_{>h}-1}\m e_{n-p_{>h}-1}+e_{n-p_{>h}},e_{n-p_{>h}-2}\m e_{n-p_{>h}-2}+e_{n-p_{>h}-1},\\&
\do,e_{n-p_{>h}-p_h+1}\m e_{n-p_{>h}-p_h+1}+e_{n-p_{>h}-p_h+2},\\&
e_0\m e_0-2e'_{n-p_{>h}}+2e'_{n-p_{>h}-1}-2e'_{n-p_{>h}-2}+\do+(-1)^{p_h}2e'_{n-p_{>h}-p_h+1},\\&
e'_{n-p_{>h}}\m e'_{n-p_{>h}}-e'_{n-p_{>h}-1}+\do+(-1)^{p_h-1}e'_{n-p_{>h}-p_h+1},\\&
e'_{n-p_{>h}-1}\m e'_{n-p_{>h}-1}-e'_{n-p_{>h}-2}+\do+(-1)^{p_h-2}e'_{n-p_{>h}-p_h+1},\\&
\do,e'_{n-p_{>h}-p_h+2}\m e'_{n-p_{>h}-p_h+2}-e'_{n-p_{>h}-p_h+1}].\endalign$$
Hence $u_w\i$ is given by
$$\align&e_i\m e_i+e_{i+1}\text{ if }i\in[1,n],i\n\{p_1,p_1+p_2,\do,p_1+p_2+\do+p_\s\},\\&
e_{p_1+p_2+\do+p_r}\m e_{p_1+p_2+\do+p_r}+e_0+\sum_{v\in[1,p_r]}(-1)^ve'_{p_1+p_2+\do+p_r-v+1}
\\&+2\sum_{v\in[1,p_{r-1}]}(-1)^ve'_{p_1+p_2+\do+p_{r-1}-v+1}+\do+2\sum_{v\in[1,p_1]}(-1)^ve'_{p_1-v+1}
\text{ if }r\in[1,\s],\\&e_0\m e_0+2\sum_{v\in[1,p_\s]}(-1)^ve'_{p_1+p_2+\do+p_\s-v+1}\\&
+2\sum_{v\in[1,p_{\s-1}]}(-1)^ve'_{p_1+p_2+\do+p_{\s-1}-v+1}+\do+2\sum_{v\in[1,p_1]}(-1)^ve'_{p_1-v+1},\\&
e'_{p_1+p_2+\do+p_r-j}\m\sum_{v\in[1,p_r-j]}(-1)^ve'_{p_1+p_2+\do+p_r-j-v+1}\text{ if }r\in[1,\s],j\in[0,p_r-2],
\\&e'_{i_1+i_2+\do+p_{r-1}+1}\m e'_{p_1+p_2+\do+p_{r-1}+1} \text{ if }r\in[1,\s].\endalign$$
We set
$$\align&e^t_j=e_{p_1+\do+p_{t-1}+j}, (t\in[1,\s],j\in[1,p_t]),\\&
e'{}^t_1=\sum_{v\in[1,p_t]}(-1)^ve'_{p_1+p_2+\do+p_t-v+1}, (t\in[1,\s]),\\&
e'{}^t_j=N^{j-1}e'{}^t_1, (t\in[1,\s],j\ge2).\endalign$$
Note that $e'{}^t_j=0$ if $j>p_t$. Clearly for any $t\in[1,\s]$, the elements $e'{}^t_j (j\in[1,p_t])$ span the 
same subspace as the elements $e'_{p_1+p_2+\do+p_t-v+1} (v\in[1,p_t])$. It follows that 
$e^t_j,e'{}^t_j (t\in[1,\s],j\in[1,p_t])$ and $e_0$ form a basis of $V$. In this basis the action of $N$ is given
by
$$\align&e^t_j\m e^t_{j+1}, (t\in[1,\s],j\in[1,p_t-1]),\\&
e^t_{p_t}\m e_0+e'{}^t_1+2e'{}^{t-1}_1+\do+2e'{}^1_1, (t\in[1,\s]),\\&e_0\m2e'{}^1_1+2e'{}^2_1+\do+2e'{}^\s_1,\\&
e'{}^t_j\m e'{}^t_{j+1}, (t\in[1,\s], j\in[1,p_t-1]),\\&e'{}^t_{p_t}\m0,(t\in[1,\s]).\endalign$$
For $t\in[1,\s],j\ge1$, we set $f^t_j=e'{}^t_j+2\sum_{t'\in[t+1,s]}e'{}^{t'}_j,(t\in[1,\s])$ and 
$\e=e_0+2\sum_{t\in[1,\s]}e'{}^t_1$. Clearly, 

(a) {\it  $e^t_j,f^t_j (t\in[1,\s],j\in[1,p_t])$ and $\e$ form a basis of $V$.}
\nl
In this basis the action of $N$ is given by
$$\align&e^t_j\m e^t_{j+1}, (t\in[1,\s],j\in[1,p_t-1]),\\&e^t_{p_t}\m\e-f^t_1, (t\in[1,\s]),\\&
\e\m2z_1+2z_2,\\&f^t_j\m f^t_{j+1}, (t\in[1,\s],j\in[1,p_t-1]),\\&f^t_{p_t}\m0\endalign$$
where for $j\ge1$ we set 
$$z_j=-\sum_{t\in[1,\s]}(-1)^tf^t_j.$$
(We use that $-\sum_{t\in[1,\s]}(-1)^tf^t_1=\sum_{t\in[1,\s]}e'{}^t_1$.) In the case where $\k_\s=0$ we set 
$$\Xi=-2\sum_{x\in[1,\s]}(-1)^xe^x_{p_x}+\e-2z_1.$$
We have
$$\align&\Xi=-2\sum_{x\in[1,\s]}(-1)^xe^x_{p_x}+e_0+2\sum_{t\in[1,\s]}e'{}^t_1
+2\sum_{t\in[1,\s]}(-1)^tf^t_1\\&=-2\sum_{x\in[1,\s]}(-1)^xe^x_{p_x}+e_0.\endalign$$
From the last expression for $\Xi$ we see that $Q(\Xi)=1$. 

\subhead 2.8\endsubhead
In the setup of 2.7 we assume that $p=2$. Then the action of $N$ in the basis 2.7(a) is given by
$$e^t_1\m e^t_2\m\do\m e^t_{p_t}\m f^t_1\m f^t_2\m\do\m f^t_{p_t}\m0\text{ if }t\in[1,\s];\e\m0.$$
Thus the Jordan blocks of $N:V@>>>V$ have sizes $2p_1,2p_2,\do,2p_\s,1$.

\subhead 2.9\endsubhead
In the setup of 2.7 we assume that $p\ne2$. For $t\in[1,\s],j\ge1$, we set $E^t_j=N^{j-1}e^t_1$. Let $\cv$ be the
subspace of $V$ spanned by the vectors $E^t_j (t\in[1,\s],j\ge1)$. Clearly, $N\cv\sub\cv$. We show:

(i) {\it  if $\k_\s=1$ then $\cv=V$;}

(ii) {\it if $\k_\s=0$ then $\cv$ is equal to $\cv'$, the codimension $1$ subspace of $V$ with basis 
$e^t_j(t\in[1,\s],j\in[1,p_t])$, $f^t_j (t\in[1,\s],j\in[2,p_t])$, $\e-f^t_1 (t\in[1,\s])$ (note that this
subspace contains $z_1$).}
\nl
For $t\in[1,\s]$ we have 

(iii) $E^t_j=e^t_j$ if $j\in[1,p_t]$, $E^t_{p_t+1}=\e-f^t_1$, $E^t_{p_t+a}=2z_{a-1}+2z_a-f^t_a$ if $a\ge2$.
\nl
Define a linear map $\ph:V@>>>\kk$ by $e^t_j\m0 (t\in[1,\s],j\in[1,p_t])$, $f^t_j\m0 (t\in[1,\s],j\in[2,p_t])$, 
$f^t_1\m1 (t\in[1,\s])$, $\e\m1$. Then $\ph(E^t_j)=0$ if $t\in[1,\s],j\in[1,p_t]$, $\ph(E^t_{p_t+1})=0$, 
$\ph(E^t_{p_t+a})=0$ if $t\in[1,\s],a\ge3$, $\ph(E^t_{p_t+2})=2\ph(z_1)=2\k_\s$ if $t\in[1,\s]$. If $\k_\s=0$, 
the last expression is $0$ so that $\cv\sub\ker\ph$ and $\cv\ne V$.

In any case from (iii) we see that $e^t_j\in\cv$ for any $t\in[1,\s],j\in[1,p_t]$ and $\e-f^t_1\in\cv$ for any 
$t\in[1,\s]$. We have $-\sum_{t\in[1,\s]}(-1)^t(\e-f^t_1)=\k_\s\e-z_1\in\cv$. For $j\ge2$ we have (using (iii)) 
$$-\sum_{t\in[1,\s]}(-1)^t(2z_{j-1}+2z_j-f^t_j)=\k_\s(2z_{j-1}+2z_j)-z_j\in\cv.$$
Thus if $\k_\s=0$ we have $z_j\in\cv$ for all $j\ge1$ hence $f^t_j\in\cv$ for all $t\in[1,\s],j\ge2$; hence 
$\cv'\sub\cv$. Since $\codim_V\cv'=1$ and $\codim_V\cv\ge2$, it follows that $\cv=\cv'$. Now assume that 
$\k_\s=1$. For $j\ge1$ we have $2z_j+z_{j+1}\in\cv$ hence 
$$z_j=2\i(2z_j+z_{j+1})-2^{-2}(2z_{j+1}+z_{j+2})+2^{-3}(2z_{j+2}+z_{j+3})-\do\in\cv.$$
It follows that $f^t_j\in\cv$ for $t\in[1,\s],j\ge2$. We have $\e-z_1\in\cv$ hence $\e\in\cv$ and $f^t_1\in\cv$
for $t\in[1,\s]$. We see that $\cv=V$. This proves (i) and (ii).

\subhead 2.10\endsubhead
In the setup of 2.7, we assume that $p\ne2$. Recall that $\ps:[1,\s]@>>>\{-1,0,1,\}$ is defined in 1.6.
For $t\in[1,\s],j\ge1$ we define $\tE^t_j\in V$ by
$$\tE^t_j=E^t_j-2\sum_{x\in[1,t-1]}(-1)^xE^x_{p_x-p_t+j}-2\sum_{x\in[1,t-1]}(-1)^xE^x_{p_x-p_t+j-1}$$
if $\ps(t)=1$ (here the last $E^x_{p_x-p_t+j-1}$ is defined since $p_x-p_t+j-1\ge0$ if $p_x>p_t,j\ge1$);
$$\tE^t_j=E^t_j-E^{t-1}_j$$   
if $t$ is odd and $\ps(t)=0$ (in this case we necessarily have $t>1$ hence $E^{t-1}_j$ is defined);
$$\tE^t_j=E^t_j+\sum_{x\in[1,t-1]}(-1)^xE^x_{p_x-p_t+j}
-\sum_{v\ge1}(-2)^{-v+1}\sum_{x\in[1,t-1]}(-1)^xE^x_{p_x-p_t+v+j-1}$$
if $t$ is even. We show:
$$\text{ if $t\in[1,\s]$, $j\ge2p_t+\ps(t)+1$, then }\tE^t_j=0.\tag a$$
Case 1. Assume that $\ps(t)=1$. In this case we have $j\ge 2p_t+2$. Hence  $j\ge p_t+3$ and for any $x\in[1,t-1]$
we have $p_x-p_t+j-1\ge p_x+2$. Thus
$$\align&\tE^t_j=(2z_{j-p_t-1}+2z_{j-p_t}-f^t_{j-p_t})
-2\sum_{x\in[1,t-1]}(-1)^x(2z_{j-p_t-1}+2z_{j-p_t}-f^x_{j-p_t})\\&
-2\sum_{x\in[1,t-1]}(-1)^x(2z_{j-p_t-2}+2z_{j-p_t-1}-f^x_{j-p_t-1})\\&=(2z_{j-p_t-1}+2z_{j-p_t}-f^t_{j-p_t})
+2\sum_{x\in[1,t-1]}(-1)^xf^x_{j-p_t}+2\sum_{x\in[1,t-1]}(-1)^xf^x_{j-p_t-1}\\&
=(2z_{j-p_t-1}+2z_{j-p_t}-f^t_{j-p_t})-2z_{j-p_t}-2\sum_{x\in[t,\s]}(-1)^xf^x_{j-p_t}\\&
-2z_{j-p_t-1}-2\sum_{x\in[t,\s]}(-1)^xf^x_{j-p_t-1}\\&
=-f^t_{j-p_t}-2\sum_{x\in[t,\s]}(-1)^xf^x_{j-p_t}-2\sum_{x\in[t,\s]}(-1)^xf^x_{j-p_t-1}.\endalign$$
This is zero since $j-p_t\ge p_t+1$ and for any $x\in[t,\s]$ we have $j-p_t-1\ge p_x+1$.

Case 2. Assume that $t$ is odd and $\ps(t)=0$. In this case we have $t>1$ and $p_t=p_{t-1}$. We have 
$j\ge 2p_t+1$. Hence $j\ge p_t+2=p_{t-1}+2$. Thus
$$\align&\tE^t_j=(2z_{j-p_t-1}+2z_{j-p_t}-f^t_{j-p_t})-(2z_{j-p_{t-1}-1}+2z_{j-p_{t-1}}-f^{t-1}_{j-p_{t-1}})\\&=
-f^t_{j-p_t}+f^{t-1}_{j-p_{t-1}}.\endalign$$
This is zero since $j-p_t\ge p_t+1=p_{t-1}+1$.

Case 3. Assume that $\ps(t)=-1$ and $p_t>1$. In this case we have $j\ge2p_t$ hence $j\ge p_t+2$ and 
$p_x-p_t+v+j-1\ge p_x+2$ (if $x\in[1,t-1],v\ge1$). Thus
$$\align&\tE^t_j=(2z_{j-p_t-1}+2z_{j-p_t}-f^t_{j-p_t})
+\sum_{x\in[1,t-1]}(-1)^x(2z_{-p_t+j-1}+2z_{-p_t+j}-f^x_{-p_t+j})\\&
-\sum_{v\ge2}(-2)^{-v+1}\sum_{x\in[1,t-1]}(-1)^x(2z_{-p_t+v+j-2}+2z_{-p_t+v+j-1}-f^x_{-p_t+v+j-1})\\&
=(2z_{j-p_t-1}+2z_{j-p_t}-f^t_{j-p_t})-(2z_{-p_t+j-1}+2z_{-p_t+j})-\sum_{x\in[1,t-1]}(-1)^xf^x_{-p_t+j}\\&
+\sum_{v\ge2}(-2)^{-v+1}(2z_{-p_t+v+j-2}+2z_{-p_t+v+j-1})\\&
+\sum_{v\ge2}(-2)^{-v+1}\sum_{x\in[1,t-1]}(-1)^xf^x_{-p_t+v+j-1}\\&=(2z_{j-p_t-1}+2z_{j-p_t}-f^t_{j-p_t})\\&
-(2z_{-p_t+j-1}+2z_{-p_t+j})+\sum_{v\ge2}(-2)^{-v+1}(2z_{-p_t+v+j-2}+2z_{-p_t+v+j-1})\\&
+z_{-p_t+j}+\sum_{x\in[t,\s]}(-1)^xf^x_{-p_t+j}\\&
-\sum_{v\ge2}(-2)^{-v+1}z_{-p_t+v+j-1}-\sum_{v\ge2}(-2)^{-v+1}\sum_{x\in[t,\s]}(-1)^xf^x_{-p_t+v+j-1}\\&
=z_{-p_t+j}+\sum_{v\ge2}(-2)^{-v+1}(2z_{-p_t+v+j-2}+z_{-p_t+v+j-1})\\&
+\sum_{x\in[t+1,\s]}(-1)^xf^x_{-p_t+j}-\sum_{v\ge2}(-2)^{-v+1}\sum_{x\in[t,\s]}(-1)^xf^x_{-p_t+v+j-1}\\&
=\sum_{x\in[t+1,\s]}(-1)^xf^x_{-p_t+j}-\sum_{v\ge2}(-2)^{-v+1}\sum_{x\in[t,\s]}(-1)^xf^x_{-p_t+v+j-1}.\endalign$$
This is zero: for $v\ge2$ and $x$ in the second sum we have $-p_t+v+j-1\ge p_x+1$ (since $j\ge 2p_t\ge p_t+p_x$);
for $x$ in the first sum we have $-p_t+j\ge p_x+1$ (since $j\ge2p_t>p_t+p_x$).

Case 4. Assume that $\ps(t)=-1$ and $p_t=1$. Since for any $x\in[t+1,\s]$ we have $p_x<p_t$ hence $p_x=0$ we see 
that $[t+1,\s]=\em$ hence $t=\s$ is even. We can assume that $j=2p_t=2$. We have
$$\align&\tE^t_j=E^\s_2+\sum_{x\in[1,\s-1]}(-1)^xE^x_{p_x+v}
-\sum_{v\ge1}(-2)^{-v+1}\sum_{x\in[1,\s-1]}(-1)^xE^x_{p_x+v}\\&=\e-f^\s_1+\sum_{x\in[1,\s-1]}(-1)^x(\e-f^x_1)
\\&-\sum_{v\ge2}(-2)^{-v+1}\sum_{x\in[1,\s-1]}(-1)^x(2z_{v-1}+2z_v-f^x_v)\\&
=z_1-\sum_{v\ge2}(-2)^{-v+1}\sum_{x\in[1,\s-1]}(-1)^x(2z_{v-1}+2z_v-f^x_v)\\&
=z_1-\sum_{v\ge2}(-2)^{-v+1}(-2z_{v-1}-2z_v)+\sum_{v\ge2}(-2)^{-v+1}z_v\\&
=z_1+\sum_{v\ge2}(-2)^{-v+1}(2z_{v-1}+z_v)=0.\endalign$$
(We have used that $f^\s_2=0$.)

Case 5. Assume that $t$ is even and $\ps(t)=0$. In this case we can have $j\ge2p_t+1$ hence $j\ge p_t+2$ and 
$p_x-p_t+v+j-1\ge p_x+2$ (if $x\in[1,t-1],v\ge1$). By the same computation as in Case 3 we have
$$\tE^t_j=\sum_{x\in[t+1,\s]}(-1)^xf^x_{-p_t+j}-\sum_{v\ge2}(-2)^{-v+1}\sum_{x\in[t,\s]}(-1)^xf^x_{-p_t+v+j-1}.$$
This is zero: for $v\ge2$ and $x$ in the second sum we have $-p_t+v+j-1\ge p_x+1$ (since 
$j\ge 2p_t+1\ge p_t+p_x$); for $x$ in the first sum we have $-p_t+j\ge p_x+1$ (since $j\ge2p_t+1\ge p_t+p_x+1$).

Note that if $\k_\s=0$ we have
$$N\Xi=-2\sum_{x\in[1,\s]}(-1)^x(\e-f^x_1)+2z_1+2z_2-2z_2=0.\tag b$$
We show:

(c) {\it if $\k_\s=1$, the elements $\tE^t_j (t\in[1,\s],j\in[1,2p_t+\ps(t)])$ form a basis of $V$; if $\k_\s=0$, 
the elements $\tE^t_j (t\in[1,\s],j\in[1,2p_t+\ps(t)])$ and $\Xi$ form a basis of $V$.}
\nl
Since $\tE^t_j$ is equal to $E^t_j$ plus a linear combination of elements $E^x_{j'}$ with $x\in[1,t-1]$ and 
$E^t_j=0$ for large $j$ we see that the subspace of V spanned by $\tE^t_j (t\in[1,\s],j\ge1)$ coincides with the 
subspace of V spanned by $E^t_j (t\in[1,\s],j\ge1)$, that is with $\cv$ (see 2.9). Using 2.9 we see that if
$\k_\s=1$, this subspace is equal to $V$ and that if $\k_\s=0$, this subspace together with $\Xi$ spans $V$. Using
this and (b) we see that if $\k_\s=1$, the elements $\tE^t_j (t\in[1,\s],j\in[1,2p_t+\ps(t)])$ span $V$; if 
$\k_\s=0$, the elements $\tE^t_j (t\in[1,\s],j\in[1,2p_t+\ps(t)])$ and $\Xi$ span $V$. It is then enough to show 
that $\sum_{t\in[1,\s]}(2p_t+\ps(t))+(1-\k_\s)=2n+1$. This follows from 1.6(b).

Using (a),(b) we see that the action of $N$ on the elements in the basis of $V$ described in (c) is as follows:
$$\tE^t_1\m\tE^t_2\m\do\m\tE^t_{2p_t+\ps(t)}\m0 (t\in[1,\s])\text{ and if $\k_\s=0$}, \Xi\m0.$$
Thus the Jordan blocks of $N:V@>>>V$ have sizes $2p_1+\ps(1),2p_2+\ps(2),\do,2p_\s+\ps(\s)$ (and $1$, if 
$\k_\s=0$). 

\subhead 2.11\endsubhead
Assume now that $\k=1,Q\ne0$, $\nn\ge5$. Let $U$ be a codimension $1$ subspace of $V$ such that $Q|_U$ is a 
nondegenerate quadratic form. Then $U$ together with the restriction of $(,)$ and $Q$ is as in 1.3 (with $U$
instead of $V$). Define $\cf_0$ in terms of $U$ in the same way as $\cf$ was defined in terms of $V$. We define a
map $\io:\cf_0@>>>\cf$ by $(0=U_0\sub U_1\sub\do U_{2n}=U)\m(0=V_0\sub V_1\sub\do V_{2n+1}=V)$ where $V_i=U_i$ for
$i\in[1,n]$ (note that $V_i$ are then uniquely determined for $i\in[n+1,2n+1]$. This is an imbedding. If 
$U_*\in\cf_0$, $U'_*\in\cf_0$ then $a_{U_*,U'_*}\in W$ is defined as in 1.4 (with $U$ instead of $V$). Let 
$V_*=\io(U_*),V'_*=\io(U'_*)$. Then $a_{V_*,V'_*}\in W$ is defined as in 1.4. From the definitions we see that
$a_{U_*,U'_*}=a_{V_*,V'_*}$.

\subhead 2.12\endsubhead
Assume that $\k=1,Q\ne0,\nn\ge5$. Let $p_*=(p_1\ge p_2\ge\do\ge p_\s)\in\cp^+_n$. Let $w\in W'$ and $u_w,N$ be as
in 2.6. Let $U$ be the subspace of $V$ spanned by $e_i,e'_i (i\in[1,n])$ if $p=2$ 
and let $U=(\kk\Xi)^\pe$ if $p\ne2$. Then $\dim U=2n$ and $Q|_U$ is a nondegenerate quadratic form with associated
bilinear form $(,)|_U$. (If $p\ne2$ we use that $Q(\Xi)=1$.) Note that $U$ is $u_w$-stable and $u_w$ is a 
unipotent isometry of $U$. (If $p\ne2$ we use that $N\Xi=0$.) For $i\in[1,n]$ we have $e_i\in U$ hence 
$\fV_i\sub U$. Now $\io:\cf_0@>>>\cf$ is defined as in 2.11. Define 
$U_*=(0=U_0\sub U_1\sub\do\sub U_{2n}=U)\in\cf_0$ where $U_i=\fV_i$ ($i\in[1,n])$. Let $U'_*=u_w\i U_*\in\cf_0$. 
Let $V_*=\io(U_*)\in\cf$, $V'_*=\io(U'_*)$. Then $V_*=(0=V_0\sub V_1\sub\do\sub V_{2n+1}=V)$ where $V_i=\fV_i$ 
for $i\in[1,n]$ and $V'_*=u_w\i(V_*)$. By 2.4(b) we have $a_{V_*,V'_*}=w_{p_*}\i\in W$. By 2.11 we have
$a_{U_*,U'_*}=a_{V_*,V'_*}$. It follows that $a_{U_*,U_*}=w_{p_*}\i\in W$. Since 
$w_{p_*}\in W'$ we see that $u_w|_U$ automatically belongs to the identity component of the isometry group of $U$.
If $p\ne2$ the Jordan blocks of $N:U@>>>U$ have sizes $2p_1+\ps(1),2p_2+\ps(2),\do,2p_\s+\ps(\s)$; indeed, by 
2.10, these are the sizes of the Jordan blocks of $N:V/\kk\Xi@>>>V/\kk\Xi$ which can be identified with $N:U@>>>U$
since the direct sum decomposition $V=U\op\kk\Xi$ is $N$-stable. If $p=2$ the Jordan blocks of $N:U@>>>U$ have 
sizes $2p_1,2p_2,\do,2p_\s$; indeed, by 2.8, these are the sizes of the Jordan blocks of 
$N:V/\kk e_0@>>>V/\kk e_0$ which can be identified with $N:U@>>>U$ since the direct sum decomposition 
$V=U\op\kk e_0$ is $N$-stable.)

\head 3. Isometry groups\endhead
\proclaim{Proposition 3.1} Let $V$ be a $\kk$-vector space of finite dimension. Let $g\in GL(V)$ be a unipotent 
element with Jordan blocks of sizes $n_1\ge n_2\ge\do\ge n_u\ge1$. We set $n_i=0$ for $i>u$. Let 
$m_1\ge m_2\ge\do\ge m_f\ge1$ be integers such that $m_1+m_2+\do+m_f=\dim V$. Assume that there exist vectors 
$x_1,x_2,\do,x_f$ in $V$ such that the vectors $g^ix_r (r\in[1,f],i\in[0,m_r-1])$ form a basis of $V$. We set 
$m_i=0$ for $i>f$. For any $c\ge1$ we have 
$$m_1+m_2+\do+m_c\le n_1+n_2+\do+n_c.\tag a$$
In particular we have $u\le f$.
\endproclaim
We set $N=g-1\in\End(V)$. For $r\in[1,f]$ and $i\in[0,m_r-1]$ we set $v_{r,i}=N^ix_r$. Note that 
$(v_{r,i})_{r\in[1,f],i\in[0,m_r-1]}$ is a basis of $V$. Note that for any $k\ge0$, the subspace $N^k(V)$ contains
the vectors $v_{r,i}(r\in[1,f],i\in[0,m_r-1],i\ge k)$. Hence 
$$\dim N^kV\ge\sum_{r\ge1}\max(m_r-k,0).\tag b$$
By assumption we can find a basis $(v'_{r,i})_{r\in[1,u],i\in[0,n_r-1]}$ of $V$ such that for any $r\in[1,u]$ we 
have $v'_{r,i}=Nv'_{r,i-1}$ if $i\in[1,n_r-1]$, $Nv'_{r,n_r-1}=0$. Then for any $k\ge0$ the subspace $N^kV$ is 
spanned by $\{v'_{r,i};r\in[1,u],i\in[0,n_r-1],i\ge k\}$. Hence $\dim N^k(V)=\sum_{r\ge1}\max(n_r-k,0)$. Thus for
any $k\ge0$ we have 
$$\sum_{r\ge1}\max(n_r-k,0)\ge\sum_{r\ge1}\max(m_r-k,0).\tag c$$
We show (a) by induction on $c$. Assume that $m_1>n_1$. The left hand side of (c) with $k=n_1$ is 
$\sum_{r\ge1}\max(n_r-n_1,0)=0$; the right hand side is 

$m_1-n_1+\sum_{r\ge2}\max(m_r-k,0)\ge m_1-n_1>0$;
\nl
thus we have $0\ge m_1-n_1>0$, a contradiction. Thus (a) holds for $c=1$. Assume now that $c\ge2$ and that (a) 
holds when $c$ is replaced by $c-1$. Assume that $m_1+\do+m_c>n_1+\do+n_c$. Then 

$m_c-n_c>(n_1+\do+n_{c-1})-(m_1+\do+m_{c-1})\ge0$.
\nl
Hence 

$\max(m_c-n_c,0)>(n_1+\do+n_{c-1})-(m_1+\do+m_{c-1})$. 
\nl
The left hand side of (c) with $k=n_c$ is 

$\sum_{r\ge1}\max(n_r-n_c,0)=\sum_{r\in[1,c-1]}(n_r-n_c)$;
\nl
the right hand side is
$$\align&\sum_{r\ge1}\max(m_r-n_c,0)\ge\sum_{r\in[1,c-1]}\max(m_r-n_c,0)+(m_c-n_c)>\\&
\sum_{r\in[1,c-1]}\max(m_r-n_c,0)+(n_1+\do+n_{c-1})-(m_1+\do+m_{c-1})\\&
\ge\sum_{r\in[1,c-1]}(m_r-n_c)+(n_1+\do+n_{c-1})-(m_1+\do+m_{c-1})\\&=(n_1+\do+n_{c-1})-(c-1)n_c.\endalign$$
Thus (c) implies $\sum_{r\in[1,c-1]}(n_r-n_c)>(n_1+\do+n_{c-1})-(c-1)n_c$. This is a contradiction. We see that 
$m_1+\do+m_c\le n_1+\do+n_c$. This yields the induction step.  This proves (a).

We have $n_1+n_2+\do+n_u=\dim V=m_1+m_2+...+m_f\le n_1+n_2+\do+n_f$. Hence $n_1+n_2+\do+n_u\le n_1+n_2+\do+n_f$. 
If $u>f$ we deduce $n_{f+1}+\do+n_u\le 0$ hence $n_{f+1}=\do=n_u=0$, absurd. Thus $u\le f$. The proposition is
proved.

\subhead 3.2 \endsubhead
In the remainder of this section $V,\nn,n,(,),Q,\cf,\k,Is(V),G$ are as in 1.3. For any sequence 
$v_1,v_2,\do,v_s$ in $V$ let $S(v_1,v_2,\do,v_s)$ be the subspace of $V$ spanned by $v_1,v_2,\do,v_s$. 

Let $p_*=(p_1\ge p_2\ge\do\ge p_\s)\in\cp_n$. Let $(V_*,V'_*)\in\cf\T\cf$ be such that $a_{V_*,V'_*}=w_{p_*}$ 
(see 1.6). From the definitions we see that for any $r\in[1,\s]$ we have
$$\dim(V'_{p_{<r}+i}\cap V_{p_{<r}+i})=p_{<r}+i-r,\qua \dim(V'_{p_{<r}+i}\cap V_{p_{<r}+i+1})=p_{<r}+i-r+1\tag a$$
if $i\in[1,p_r-1]$;
$$\dim(V'_{p_{\le r}}\cap V_{\nn-p_{<r}-1})=p_{\le r}-r,\qua \dim(V'_{p_{\le r}}\cap V_{\nn-p_{<r}})
=p_{\le r}-r+1.\tag b$$
In this setup we have the following result.

\proclaim{Proposition 3.3} Let $g\in Is(V)$ be such that $gV_*=V'_*$. There exist vectors $v_1,v_2,\do,v_\s$ in 
$V$ (each $v_i$ being unique up to multiplication by $\pm1$) such that, setting 
$Z_k=S(v_k,gv_k,\do,g^{p_k-1}v_k)$ for $k\in[1,\s]$, the following hold for any $r\in[1,\s]$:

(i) $V_{p_{<r}+i}=Z_1+\do+Z_{r-1}+S(v_r,gv_r,\do,g^{i-1}v_r)$ for $i\in[0,p_r]$;

(ii) $(g^iv_t,v_r)=0$ for any $1\le t<r$, $i\in[-p_t,p_t-1]$;

(iii) $(v_r,g^iv_r)=0$ for $i\in[-p_r+1,p_r-1]$, $Q(v_r)=0$ and $(v_r,g^{p_r}v_r)=1$;

(iv) the vectors $(g^{-p_t+i}v_t)_{t\in[1,r],i\in[0,2p_t-1]}$ are linearly independent;

(v) setting $E_r=S(g^{-p_t+i}v_t;t\in[1,r],i\in[0,p_t-1])$ we have $V=V_{p_{\le r}}\op E_r^\pe$.
\nl
If $\k=1$ there exists a vector $v_{\s+1}\in V$ (unique up to multiplication by $\pm1$) such that

(ii${}'$) $(g^iv_t,v_{\s+1})=0$ for any $1\le t<s+1$, $i\in[-p_t,p_t-1]$;

(iii${}'$) $Q(v_{\s+1})=1$.
\nl
Moreover,

(vi) if $\k=0$, the vectors $(g^jv_t)_{t\in[1,\s],j\in[-p_t,p_t-1]}$ form a basis of $V$; if $\k=1$, the vectors 
$(g^jv_t)_{t\in[1,\s],j\in[-p_t,p_t-1]}$ together with $v_{\s+1}$ form a basis of $V$.
\endproclaim
We shall prove the following statement for $u\in[1,\s]$.

(a) {\it There exist vectors $v_1,v_2,\do,v_u$ in $V$ (each $v_i$ being unique up to multiplication by $\pm1$) 
such that for any $r\in[1,u]$, (i)-(v) hold.}
\nl
We can assume that (a) holds when $u$ is replaced by a strictly smaller number in $[1,u]$. (This assumption is
empty when $u=1$.) In particular $v_1,\do,v_{u-1}$ are defined. By assumption we have 
$V=V_{p_{<u}}\op E_{u-1}^\pe$ hence $V_{p_{<u}+1}\cap E_{u-1}^\pe$ is a line. (We set $E_0=0$ so that 
$E_0^\pe=V$.) Let $v_u$ be a nonzero vector on this line. 

We show that (ii) holds for $r\in[1,u]$. It is enough to show this when $r=u$. From the definition we have 
$V_{p_{<u}+1}=V_{p_{<u}}+\kk v_u=Z_1+\do+Z_{u-1}+\kk v_u$. Since $p_{<u}+1\le\nn/2$, $V_{p_{<u}+1}$ is isotropic.
Hence for $t\in[1,u-1]$ we have $(Z_t,v_u)=0$ that is $(g^iv_t,v_u)=0$ for any $i\in[0,p_t-1]$. From the 
definition of $v_u$ we have $(g^iv_t,v_u)=0$ for any $t\in[1,u-1]$ and any $i\in[-p_t,-1]$. Thus, (ii) holds when
$r=u$.

We show that (i) holds for $r\in[1,u]$. It is enough to show (i) when $r=u$ and $i\in[0,p_u]$. We argue by 
induction on $i$. For $i=0$ the result follows from the induction hypothesis. For $i=1$ the result follows from 
the definition of $v_u$. Assume now that $i\in[2,p_u]$. Let $j=i-1$. By the induction hypothesis we have 
$V_{p_{<u}+j}=Z_1+\do+Z_{u-1}+S(v_u,gv_u,\do,g^{j-1}v_u)$, hence 
$gV_{p_{<u}+j}=gZ_1+\do+gZ_{u-1}+S(gv_u,g^2v_u,\do,g^jv_u)$ and the $p_{<u}+j$ vectors 
$$v_1,gv_1,\do,g^{p_1-1}v_1,\do,v_{u-1},gv_{u-1},\do,g^{p_{u-1}-1}v_{u-1},v_u,gv_u,\do,g^{j-1}v_u$$
form a basis of $V_{p_{<u}+j}$. Hence the $p_{<u}+j-u$ vectors
$$gv_1,\do,g^{p_1-1}v_1,\do,gv_{u-1},\do,g^{p_{u-1}-1}v_{u-1},gv_u,\do,g^{j-1}v_u$$
are linearly independent; they are contained in $gV_{p_{<u}+j}\cap V_{p_{<u}+j}$ (of dimension $p_{<u}+j-u$) hence
$$\align&S(gv_1,\do,g^{p_1-1}v_1,\do,gv_{u-1},\do,g^{p_{u-1}-1}v_{u-1},gv_u,\do,g^{j-1}v_u)\\&=
gV_{p_{<u}+j}\cap V_{p_{<u}+j}.\endalign$$
Since $\dim(gV_{p_{<u}+j}\cap V_{p_{<u}+j+1})=p_{<u}+j-u+1$, we see that
$$\align&gV_{p_{<u}+j}\cap V_{p_{<u}+j+1}\\&=
S(gv_1,\do,g^{p_1-1}v_1,\do,gv_{u-1},\do,g^{p_{u-1}-1}v_{u-1},gv_u,\do,g^{j-1}v_u,x)\endalign$$
for a unique (up to scalar) $x\in gV_{p_{<u}+j}$ of the form $x=\sum_{r=1}^{u-1}a_rg^{p_r}v_r+a_ug^jv_u$ 
where $a_r,a_u\in\kk$ are not all $0$. 

Assume that $a_r\ne0$ for some $r\in[1,u-1]$; let $r_0$ be the smallest such $r$. Note that $V_{p_{<u}+j+1}$ is an
isotropic subspace. Since $x\n V_{p_{<u}+j+1}$ and $v_{r_0}\in V_{p_{<u}+j+1}$, we have $(x,v_{r_0})=0$ hence
$\sum_{r=r_0}^{u-1}a_r(g^{p_r}v_r,v_{r_0})+a_u(g^jv_u,v_{r_0})=0$. By the definition of $v_u$ we have
$(g^jv_u,v_{r_0})=(v_u,g^{-j}v_{r_0})=0$ since $j\in[0,p_u-1]\sub[0,p_{r_0}-1]$. If $r\in[r_0+1,u-1]$ we have 
$(g^{p_r}v_r,v_{r_0})=(v_r,g^{-p_r}v_{r_0})=0$ since $-p_r\in[-p_{r_0},0]$ (we use (ii)). We see that 
$a_{r_0}(g^{p_{r_0}}v_{r_0},v_{r_0})=0$. Using that $(g^{p_{r_0}}v_{r_0},v_{r_0})\ne0$ (induction hypothesis) we
deduce $a_{r_0}=0$, a contradiction. 

We see that $a_r=0$ for any $r\in[1,u-1]$ hence $a_u\ne0$; we can assume that $a_u=1$. Thus $x=g^jv_u$ and
$$\align&gV_{p_{<u}+j}\cap V_{p_{<u}+j+1}\\&=
S(gv_1,\do,g^{p_1-1}v_1,\do,gv_{u-1},\do,g^{p_{u-1}-1}v_{u-1},gv_u,\do,g^{j-1}v_u,g^jv_u).\endalign$$
Now $gV_{p_{<u}+j}\cap V_{p_{<u}+j+1}\not\sub V_{p_{<u}+j}$ (otherwise we would have 
$gV_{p_{<u}+j}\cap V_{p_{<u}+j+1}\sub gV_{p_{<u}+j}\cap V_{p_{<u}+j}$ and passing to dimensions: 
$p_{<u}+j-u+1\le p_{<u}+j-u$, a contradiction.) Since $V_{p_{<u}+j}$ is a hyperplane in $V_{p_{<u}+j+1}$ not 
containing the subspace $gV_{p_{<u}+j}\cap V_{p_{<u}+j+1}$ of $V_{p_{<u}+j+1}$, we have 
$V_{p_{<u}+j}+(gV_{p_{<u}+j}\cap V_{p_{<u}+j+1})=V_{p_{<u}+j+1}$. It follows that
$$\align&V_{p_{<u}+j+1}=S(v_1,\do,g^{p_1-1}v_1,\do,v_{u-1},\do,g^{p_{u-1}-1}v_{u-1},v_u,\do,g^{j-1}v_u)\\&+
S(gv_1,\do,g^{p_1-1}v_1,\do,gv_{u-1},\do,g^{p_{u-1}-1}v_{u-1},gv_u,\do,g^{j-1}v_u,g^jv_u)\\&
=S(v_1,\do,g^{p_1-1}v_1,\do,v_{u-1},\do,g^{p_{u-1}-1}v_{u-1},v_u,\do,g^{j-1}v_u,g^jv_u).\endalign$$
Thus (i) holds when $r=u$.

We show that (iii) holds for $r\in[1,u]$. It is enough to show this when $r=u$. Using that 
$v_u,gv_u,\do,g^{p_u-1}v_u$ are contained in $V_{p_{\le u}}$ (see (i)) which is an isotropic subspace (since 
$p_{\le u}\le\nn/2$) we see that for $i\in[0,p_u-1]$ we have $(v_u,g^iv_u)=0$ hence also $(v_u,g^{-i}v_u)=0$; 
moreover $Q(v_u)=0$. We have $\dim(gV_{p_{\le u}}\cap V_{p_{<u}+1}^\pe)=p_{\le u}-u$ and
$gV_{p_{\le u}}=S(gv_1,\do,g^{p_1}v_1,\do,gv_u,\do,g^{p_u}v_u)$. Moreover, using the part of (iii) that is 
already known, we see that 
$$S(gv_1,\do,g^{p_1-1}v_1,\do,gv_u,\do,g^{p_u-1}v_u)$$
(of dimension $p_{\le u}-u$) is contained in $gV_{p_{\le u}}\cap V_{p_{<u}+1}^\pe$ hence is equal to 
$gV_{p_{\le u}}\cap V_{p_{<u}+1}^\pe$. Hence $g^{p_u}v_u\n V_{p_{<u}+1}^\pe$ that is 
$$g^{p_u}v_u\n S(v,gv_1,\do,g^{p_1-1}v_1,\do,v_{u-1},\do,g^{p_{u-1}}v_{u-1},v_u)^\pe.$$
Since $g^{p_u}v_u\in S(v,gv_1,\do,g^{p_1-1}v_1,\do,v_{u-1},\do,g^{p_{u-1}}v_{u-1})^\pe$ it follows that\lb
$(g^{p_u}v_u,v_u)\ne0$. Replacing $v_u$ by a scalar multiple we can assume that $(v_u,g^{p_u}v_u)=1$. Thus (iii) 
holds when $r=u$.

We show that (iv) holds for $r\in[1,u]$. It is enough to show this when $r=u$. Assume that 
$f=\sum_{r\in[1,u]}\sum_{i\in[0,2p_r-1]}c_{r,i}g^{-p_r+i}v_r$ is equal to $0$ where $c_{r,i}\in\kk$ are not all 
zero. Let $i_0=\min\{i;c_{r,i}\ne0\text{ for some }r\in[1,u]\}$. Let $X=\{r\in[1,u];c_{r,i_0}\ne0\}$. We have 
$X\ne\em$ and 
$$f=\sum_{r\in X}c_{r,i_0}g^{-p_r+i_0}v_r+\sum_{r\in[1,u]}\sum_{i\in[i_0+1,2p_r-1]}c_{r,i}g^{-p_r+i}v_r.$$
Let $r_0$ be the largest number in $X$. We have 
$$\align&0=(f,g^{i_0}v_{r_0})=\sum_{r\in X}c_{r,i_0}(g^{-p_r+i_0}v_r,g^{i_0}v_{r_0})\\&
+\sum_{r\in[1,u]}\sum_{i\in[i_0+1,2p_r-1]}c_{r,i}(g^{-p_r+i}v_r,g^{i_0}v_{r_0}).\endalign$$
If $r\in X$, $r\ne r_0$ we have $(g^{-p_r+i_0}v_r,g^{i_0}v_{r_0})=0$ (using (ii) and $r<r_0$). If $r\in[1,u]$ and
$i\in[i_0+1,2p_r-1]$ we have $(g^{-p_r+i}v_r,g^{i_0}v_{r_0})=0$ (we use (ii),(iii); note that if $r<r_0$, we have 
$-p_r+i-i_0\in[-p_r,p_r-1]$; if $r\ge r_0$ we have $p_r+i_0-i\in[-p_{r_0}+1,p_{r_0}-1]$). We see that 
$$0=c_{r_0,i_0}(g^{-p_{r_0}+i_0}v_{r_0},g^{i_0}v_{r_0})=c_{r_0,i_0}(g^{-p_{r_0}}v_{r_0},v_{r_0}).$$
Using (iii) we have $(g^{-p_{r_0}}v_{r_0},v_{r_0})\ne0$ hence $c_{r_0,i_0}=0$, a contradiction. Thus (iv) holds 
when $r=u$.

We show that (v) holds for $r\in[1,u]$. It is enough to show this when $r=u$. From (iv) with $r=u$ we see that 
$\dim(E_u)=p_{\le u}$. By (i) with $r=u$ we have $E_u\sub g^{-p_t}V_{p_{\le u}}$. Since $V_{p_{\le u}}$ is 
isotropic we see that $E_u$ is isotropic hence $\dim(E_u^\pe)=\nn-\dim(E_u)=\nn-p_{\le u}$. 
Thus $\dim(V_{p_{\le u}})+\dim(E_u^\pe)=\dim V$. It is enough to show that $V_{p_{\le u}}\cap E_u^\pe=0$. Assume 
that $f=\sum_{r\in[1,u],i\in[1,p_r]}c_{r,i}g^{p_r-i}v_r$ belongs to $E_u^\pe$ and is nonzero. Here $c_{r,i}\in\kk$
are not all zero. Let $i_0=\min\{i;c_{r,i}\ne0\text{ for some }r\in[1,u]\}$. Let $X'=\{r\in[1,u];c_{r,i_0}\ne0\}$.
We have $X'\ne\em$ and 
$$f=\sum_{r\in X'}c_{r,i_0}g^{p_r-i_0}v_r+\sum_{r\in[1,u]}\sum_{i\in[i_0+1,p_r]}c_{r,i}g^{p_r-i}v_r.$$
Let $r_0$ be the smallest number in $X'$. We have 
$$\align&0=(f,g^{-i_0}v_{r_0})=\sum_{r\in X'}c_{r,i_0}(g^{p_r-i_0}v_r,g^{-i_0}v_{r_0})\\&
+\sum_{r\in[1,u]}\sum_{i\in[i_0+1,p_r]}c_{r,i}(g^{p_r-i}v_r,g^{-i_0}v_{r_0}).\endalign$$
If $r\in X'$, $r\ne r_0$, we have $(g^{p_r-i_0}v_r,g^{-i_0}v_{r_0})=(v_r,g^{-p_r}v_{r_0})=0$ (we use (ii); note 
that $r\ge r_0$ hence $p_r\le p_{r_0}$). If $r\in[1,u]$ and $i\in[i_0+1,p_r]$, we have 
$(g^{p_r-i}v_r,g^{-i_0}v_{r_0})=0$ (we use (ii),(iii); note that if $r>r_0$ we have $-p_r+i-i_0\in[-p_{r_0},0]$; 
if $r\le r_0$ we have $p_r-i+i_0\in[0,p_r-1]$.) Thus (v) holds when $r=u$. This completes the inductive proof of 
(a). Taking $u=\s$ in (a) we obtain (i)-(v).

By (v) we have $V=V_{p_{\le\s}}\op E_\s^\pe$; hence $V_{p_{\le\s}+1}\cap E_\s^\pe$ is a line. Let $v_{\s+1}$ be a 
nonzero vector on this line. From the definition we have
$$V_{p_{\le\s}+1}=V_{p_{\le\s}}+S(v_{\s+1}).\tag b$$
Since $V_{p_{\le\s}}$ is a maximal isotropic subspace of $V$ and $(V_{p_{\le\s}},v_{\s+1})=0$ (by (b)) we see that
$Q(v_{\s+1})\ne0$. Replacing $v_{\s+1}$ by a scalar multiple we can assume that (iii${}'$) holds. Using 
$v_{\s+1}\in V_{p_{\le\s}+1}=V_{\p_{\le\s}}^\pe$ and (i) we see that $(g^iv_t,v_{\s+1})=0$ for $t\in[1,\s]$, 
$i\in[0,p_t-1]$. From $v_{\s+1}\in E_\s^\pe$ we have $(g^iv_t,v_{\s+1})=0$ for $t\in[1,\s]$, $i\in[-p_t,-1]$. Thus
(ii${}'$) holds.

If $\k=0$, (vi) follows from (iv) with $r=\s$. In the rest of the proof we assume that $\k=1$. If $p=2$ we denote
by $\o$ the unique vector in $V^\pe$ such that $Q(\o)=1$. Since $V_{p_{\le\s}+1}=V_{p_{\le\s}}^\pe$ we have 
$\o\in V_{p_{\le\s}+1}$. Clearly, $\o\in E_\s^\pe$. Hence $\o\in V_{p_{\le\s}+1}\cap E_\s^\pe$. Thus we have 
$v_{\s+1}=\o$. Returning to a general $p$ we show that (vi) holds when $\k=1$. Assume that
$$f=\sum_{r\in[1,\s],i\in[0,2p_r-1]}c_{r,i}g^{-p_r+i}v_r+c_{\s+1,0}v_{\s+1}$$
is equal to $0$ where $c_{r,i}\in\kk$ are not all zero. If $c_{\s+1,0}=0$ then we have a contradiction by (iv).
So we can assume that $c_{\s+1,0}\ne0$ or even that $c_{\s+1,0}=1$. We have 
$$0=(f,v_{\s+1})=\sum_{r\in[1,\s],i\in[0,2p_r-1]}c_{r,i}(g^{-p_r+i}v_r,v_{\s+1})+(v_{\s+1},v_{\s+1}).$$
For $r,i$ in the sum we have $-p_r\le-p_r+i\le p_r-1$ hence $(g^{-p_r+i}v_r,v_{\s+1})=0$ (see (ii${}'$)) hence 
$(v_{\s+1},v_{\s+1})=0$. If $p\ne2$ we have $Q(v_{\s+1})\ne0$ hence $(v_{\s+1},v_{\s+1})\ne0$, contradiction. Hence
we may assume that $p=2$ so that $v_{\s+1}=\o$. The following proof is almost a repetition of that of (iv). We have
$\sum_{r\in[1,\s],i\in[0,2p_r-1]}c_{r,i}g^{-p_r+i}v_r+\o=0$. Assume that $c_{r,i}\in\kk$ are not all zero. Let 
$i_0=\min\{i;c_{r,i}\ne0\text{ for some }r\in[1,\s]\}$. Let $X=\{r\in[1,\s];c_{r,i_0}\ne0\}$. We have $X\ne\em$ and 
$$\sum_{r\in X}c_{r,i_0}g^{-p_r+i_0}v_r+\sum_{r\in[1,\s],i\in[i_0+1,2p_r-1]}c_{r,i}g^{-p_r+i}v_r+\o=0.$$
Let $r_0$ be the largest number in $X$. We have 
$$\align&0=(0,g^{i_0}v_{r_0})=\sum_{r\in X}c_{r,i_0}(g^{-p_r+i_0}v_r,g^{i_0}v_{r_0})\\&
+\sum_{r\in[1,\s],i\in[i_0+1,2p_r-1]}c_{r,i}(g^{-p_r+i}v_r,g^{i_0}v_{r_0}).\endalign$$
If $r\in X$, $r\ne r_0$ we have $(g^{-p_r+i_0}v_r,g^{i_0}v_{r_0})=0$ (using (ii) and $r<r_0$). If $r\in[1,\s]$ 
and $i\in[i_0+1,2p_r-1]$ we have $(g^{-p_r+i}v_r,g^{i_0}v_{r_0})=0$ (we use (ii),(iii); note that if 
$r<r_0$, we have $-p_r+i-i_0\in[-p_r,p_r-1]$; if $r\ge r_0$ we have $p_r+i_0-i\in[-p_{r_0}+1,p_{r_0}-1]$). We see
that 
$$0=c_{r_0,i_0}(g^{-p_{r_0}+i_0}v_{r_0},g^{i_0}v_{r_0})=c_{r_0,i_0}$$
(we have used (iii)) hence $c_{r_0,i_0}=0$, a contradiction. We see that $c_{r,i}=0$ for all 
$r\in[1,\s],i\in[0,2p_r-1]$. Hence $\o=0$ contradiction. This proves (vi). The proposition is proved.

\subhead 3.4\endsubhead
We preserve the setup of 3.3. For any $r\in[1,\s]$ let $X_r$ be the subspace of $V$ spanned by 
$(g^{-p_r+i}v_r)_{i\in[0,2p_r-1]}$. Let $X_{\s+1}$ be $0$ (if $\k=0$) and the subspace spanned by $v_{\s+1}$ (if 
$\k=1$). From 3.3(vi) we see that
$$V=\op_{r\in[1,\s+1]}X_r.\tag a$$
For $t\in[1,\s]$ let $W_t=\op_{r\in[1,t]}X_r$, $W'_t=\op_{r\in[t+1,\s+1]}X_r$. From (a) we see that 
$$V=W_t\op W'_t.\tag b$$
We show:
$$V=X_r\op X_r^\pe\text{ if }r\in[1,\s].\tag c$$
Let $f=\sum_{i\in[0,2p_r-1]}c_{r,i}g^{-p_r+i}v_r$ be such that $(f,g^jv_r)=0$ for any $j\in[-p_r,p_r-1]$. Here 
$c_{r,i}\in\kk$. We show that $c_{r,i}=0$ for all $i$. Assume that $c_{r,i}\ne0$ for some $i\in[0,2p_r-1]$ and let
$i_0$ be the smallest $i$ such that $c_{r,i}\ne0$. Assume first that $i_0\in[0,p_r-1]$. We have 
$$\align&0=(f,g^{i_0}v_r)=\sum_{i\in[i_0,2p_r-1]}c_{r,i}(g^{-p_r+i}v_r,g^{i_0}v_r)\\&=
c_{r,i_0}(g^{-p_r}v_r,v_r)+\sum_{i\in[i_0+1,2p_r-1]}c_{r,i}(g^{-p_r+i-i_0}v_r,v_r).\endalign$$
In the last sum we have $-p_r+1\le-p_r+i-i_0\le p_r-1$ hence the last sum is zero (see 3.3(iii)). We see that 
$c_{r,i_0}(g^{-p_r}v_r,v_r)=0$ hence $c_{r,i_0}=0$ (see 3.3(iii)), a contradiction. Thus we have 
$i_0\in[p_r,2p_r-1]$ so that $f=\sum_{i\in[p_r,2p_r-1]}c_{r,i}g^{-p_r+i}v_r$. Let $i_1$ be largest $i$ such that 
$c_{r,i}\ne0$. We have $i_1\in[p_r,2p_r-1]$ hence $-2p_r+i_1\in[-p_r,-1]$. We have 
$$\align&0=(f,g^{-2p_r+i_1}v_r)=\sum_{i\in[p_r,i_1]}c_{r,i}(g^{-p_r+i}v_r,g^{-2p_r+i_1}v_r)\\&
=c_{r,i_1}(g^{p_r}v_r,v_r)+\sum_{i\in[p_r,i_1-1]}c_{r,i}(g^{p_r+i-i_1}v_r,v_r).\endalign$$
In the last sum we have $-p_r+1\le p_r+i-i_1\le p_r-1$ hence the last sum is zero (see 3.3(iii)). We see that 
$c_{r,i_1}(g^{p_r}v_r,v_r)=0$ hence $c_{r,i_1}=0$, a contradiction. We see that $X_r\cap X_r^\pe=0$. We have 
$\dim X_r^\pe+\dim X_r\ge\dim V$ hence $V=X_r\op X_r^\pe$, as required.

\subhead 3.5\endsubhead
In the setup of 3.3 we assume that $g$ is unipotent; we set $N=g-1\in\End(V)$. We set $p_{\s+1}=\k/2$. Note that
$\p_1\ge\p_2\ge\do\ge\p_\s\ge p_{\s+1}$. For any $k\ge0$ we set
$$\L_k=\sum_{r\in[1,\s+1]}\max(2p_r-k,0).$$
We show:

(a) {\it For any $k\in\NN$ we have $\dim N^kV\ge\L_k$. Moreover, $\dim N^0V=\nn=\L_0$.}
\nl
The inequality in (a) follows from 3.1(b) using 3.3(vi). (We apply 3.1(b) with $x_1,\do,x_f$ given by 
$g^{-p_1}v_1,\do,g^{-p_\s}v_\s$, if $\k=0$, or by $g^{-p_1}v_1,\do,g^{-p_\s}v_\s,v_{\s+1}$, if $\k=1$.) The 
equality in (a) follows from $\L_0=\sum_{r\in[1,\s+1]}2p_r=\nn$. 

We now assume that $k>0$ and $d\in[1,\s]$ is such that $2p_d\ge k\ge2p_{d+1}$. Then
$\L_k=\sum_{r\in[1,d]}(2p_r-k)$. We show:

(b) {\it If $\dim N^kV=\L_k$ then $W_d,W'_d$ are $g$-stable, $W'_d=W_d^\pe$, $g:W_d@>>>W_d$ has exactly $d$ Jordan
blocks (each one has size $\ge k$) and $N^kW'_d=0$.}
\nl
For $r\in[1,\s]$ let $v'_r=g^{-p_r}v_r$; then $(g^iv'_r)_{i\in[0,2p_r-1]}$ is a basis of $X_r$ hence 
$(N^iv'_r)_{i\in[0,2p_r-1]}$ is a basis of $X_r$. For $r\in[1,d]$ let $Y_r$ be the subspace spanned by 
$N^iv'_r (i\in[k,2p_r-1])$. Note that $Y_r\sub N^kX_r$. Hence $\op_{r\in[1,d]}Y_r\sub N^kW_d\sub N^kV$. We have 
$\dim\op_{r\in[1,d]}Y_r=\sum_{r\in[1,d]}(2p_r-k)=\L_k=\dim N^kV$. Hence $\op_{r\in[1,d]}Y_r=N^kW_d=N^kV$. We have
$\op_{r\in[1,d]}Y_r\sub W_d$. Hence $N^kV\sub W_d$. We show that $NW_d\sub W_d$. Clearly $N$ maps the basis 
elements $N^iv'_r$ $(r\in[1,d],i\in[0,2p_r-2])$ into $W_d$. So it is enough to show that $N$ maps $N^{2p_r-1}v'_r$
($r\in[1,d]$) into $W_d$. But $NN^{2p_r-1}v'_r=N^{2p_r}v'_r=N^kN^{2p_r-k}v'_r\sub N^kV\sub W_d$. Thus 
$NW_d\sub W_d$. Hence $gW_d=W_d$ and $gW_d^\pe=W_d^\pe$. For $r\in[d+1,\s]$ we have $v_r\in W_d^\pe$ by 3.3(ii).
Since $gW_d^\pe=W_d^\pe$ we have $g^jv_r\in W_d^\pe$ for all $j\in\ZZ$; hence $X_r\sub W_d^\pe$. Similarly, if 
$\k=1$ we have $v_{\s+1}\in W_d^\pe$ by 3.3(ii${}'$); hence $X_{\s+1}\sub W_d^\pe$. We see that in any case 
$W'_d\sub W_d^\pe$. Since $W_d\op W'_d=V$, it follows that $W_d+W_d^\pe=V$. Since $V^\pe\sub W'_d$ we have 
$V^\pe\cap W_d=0$ hence $\dim W_d^\pe=\dim V-\dim W_d$ which, together with $W_d+W_d^\pe=V$, implies 
$W_d\op W_d^\pe=V$ and $W_d^\pe=W'_d$. In particular, $W'_d$ is $g$-stable. Let $\d$ be the number of Jordan 
blocks of $N:W_d@>>>W_d$ that is $\d=\dim(\ker N:W_d@>>>W_d)$. We have 
$\dim W_d-\d=\dim NW_d\ge\sum_{r\in[1,d]}(2p_r-1)=\dim W_d-d$. (The inequality follows from 3.1(b) applied to
$N:W_d@>>>W_d$.) Hence $\d\le d$. From the definition of $\d$ we see that $\dim(\ker N^k:W_d@>>>W_d)\le\d k$.
Recall that $\dim N^kW_d=\sum_{r\in[1,d]}(2p_r-k)=\dim W_d-kd$. Hence 
$\dim(\ker N^k:W_d@>>>W_d)=\dim W_d-\dim N^kW_d=kd$. Hence $kd\le\d k$. Since $k>0$ we deduce $d\le\d$. Hence 
$d=\d$. Since $\dim(\ker N^k:W_d@>>>W_d)=kd$ we see that each of the $d$ Jordan blocks of $N:W_d@>>>W_d$ has size
$\ge k$. Since $N^kW=N^kV$ and $V=W_d\op W'_d$ we see that $N^kW'_d=0$. Hence each Jordan block of
$N:W'_d@>>>W'_d$ has size $\le k$. This proves (b).

We show:

(c) {\it if $\dim N^kV=\L_k$ for all $k\ge0$ then for any $r\in[1,\s+1]$, $X_r$ is a $g$-stable subspace of $V$ 
and for any $r\ne r'$ in $[1,\s+1]$ we have $(X_r,X_{r'})=0$.}
\nl
Applying (b) with $k=2p_d$ for $d=1,2,\do,\s$ we see that each of the subspaces 
$X_1\sub X_1\op X_2\sub\do\sub X_1\op X_2\op\do\op X_\s$ of $V$ is $g$-stable and each of the subspaces 
$X_2\op\do\op X_{\s+1}\sps\do\sps X_\s\op X_{\s+1}\sps X_{\s+1}$ of $V$ is $g$-stable. Taking intersections we see
that each of the subspaces $X_1,X_2,\do,X_{\s+1}$ of $V$ is $g$-stable. The second assertion of (c) also follows
from (b).

\subhead 3.6\endsubhead
We preserve the setup of 3.5 and we assume that $p\ne2,Q\ne0$. For any $k>0$ such that $\dim N^kV=\L_k$ we show:

(a) {\it  if $d\in[1,\s]$ is such that $k\in[2p_{d+1},2p_d]$ then $d$ is even;}

(b) {\it $k\ne2p_r$ for $r\in[1,\s]$.}
\nl
The proof is based on the following known property of a unipotent isometry $T:W@>>>W$ of a finite dimensional 
$\kk$-vector space $W$ with a nongenerate symmetric bilinear form (assuming $p\ne2$): the number of Jordan blocks
of $T$ is congruent $\mod 2$ to $\dim W$.

By 3.5(b), $W_d$ is $g$-stable and $(,)$ is nondegenerate on $W_d$. Hence $g:W_d@>>>W_d$ has an even number of 
Jordan blocks. (Clearly, $\dim W_d$ is even.) By 3.5(b), $g:W_d@>>>W_d$ has exactly $d$ Jordan blocks. Hence $d$ 
is even, proving (a).

Assume now that $k=2p_r$ for some $r\in[1,\s]$. If $r$ is odd we have $k\in[2p_{r+1},2p_r]$ hence by (a), $r$ is 
even, a contradiction. If $r$ is even we have $r\ge2$ and $k\in[2p_r,2p_{r-1}]$ hence by (a), $r-1$ is even, a 
contradiction. This proves (b).

For any $k\ge0$ we define $\L'_k\in\NN$ by $\L'_k=\L_k+1$ if $k>0,k\in[2p_{d+1},2p_d]$ for some odd $d\in[1,\s]$ 
and $\L'_k=\L_k$ otherwise. In particular, if $k=2p_r$ for some $r\in[1,\s]$ then $\L'_k=\L_k+1$. (If $r$ is odd 
we have $k\in[2p_{r+1},p_r]$ hence $\L'_k=\L_k+1$. If $r$ is even we have $r\ge2$ and $k\in[2p_r,2p_{r-1}]$ hence 
$\L'_k=\L_k+1$.) From (a) and 3.5(a) we see:

(c) {\it  $\dim N^kV\ge\L'_k$ for any $k\ge0$; moreover $\dim N^0V=\nn=\L'_0$.}
\nl
For $r\in[1,\s+1]$ we set $\p_r=2p_r+\ps(r)$ where 

$\ps(r)=1$ if $r$ is odd, $r\le\s$ and $p_{r-1}>p_r$ (convention: $p_0=\iy$);

$\ps(r)=-1$ if $r$ is even and $p_r>p_{r+1}$ (convention: $p_{\s+2}=0$);

$\ps(r)=0$ for all other $r$.
\nl
When $r\in[1,\s]$ this definition of $\ps(r)$ agrees with that in 1.6. Note that $\ps(\s+1)$ equals $-1$ if 
$\k_\s=\k=1$ and equals $0$ otherwise.

In the remainder of this subsection we assume that 

(d) if $\k=0$ then $\k_\s=0$.
\nl
For any $k\ge0$ we set $\L''_k=\sum_{r\in[1,\s+1]}\max(\p_r-k,0)$. We show:
$$\L''_k=\L'_k\text{ for all }k\ge0.\tag e$$
Assume first that $2p_d>k>2p_{d+1}$ for some $d\in[1,\s]$. Then the conditions
$2p_r\ge k,2p_r>k,\p_r\ge k,r\le d$ are equivalent hence
$$\L''_k-\L_k=\sum_{r\in[1,\s+1];2p_r\ge k}(\p_r-k)-\sum_{r\in[1,\s+1];2p_r\ge k}(2p_r-k)=\sum_{r\in[1,d]}\ps(r).
$$
This equals $1=\L'_k-\L_k$ if $d$ is odd and equals $1+\ps(d)=1-1=0=\L'_k-\L_k$ if $d$ is even.

Next we assume that $k=2p_{d'}$ for some $d'\in[1,\s]$. There is a unique $d\in[1,\s]$ such that
$2p_{d'}=2p_d>2p_{d+1}$. The condition that $\p_r\ge k$ is equivalent to $2p_r\ge 2p_d$ (if $\ps(r)\in\{0,1\}$)
and to $2p_r\ge 2p_d$, $p_r\ne p_d$ (if $\ps(r)=-1$). Moreover for $r$ such that $p_r=p_d$, $r\ne d$ we have
$p_r=p_{r+1}$ hence $\ps(r)\ne-1$. Thus 
$$\align&\L''_k-\L_k=\sum_{r\in[1,\s+1];2p_r\ge2p_d}(\p_r-2p_d)-\sum_{r\in[1,\s+1];\ps(r)=-1;2p_r=2p_d}
(2p_r-1-2p_d)\\&-\sum_{r\in[1,\s+1];2p_r\ge2p_d}(2p_r-2p_d)
=\sum_{r\in[1,d]}\ps(r)+\sum_{r\in[1,\s+1];\ps(r)=-1;r=d}1.\endalign$$
By 1.6 this equals $(1+\ps(d))+1=(1-1)+1=1=\L'_k-\L_k$ if $d$ is even (so that $\ps(d)=-1$) and equals 
$1+0=\L'_k-\L_k$ if $d$ is odd (so that $\ps(d)\ne-1$).

Next we assume that $k\le 2p_{\s+1}$. Then $k=0$ and $\p_r\ge k,2p_r\ge k$ for all $r\in[0,1+\s]$ hence
$$\L''_k-\L_k=\sum_{r\in[1,\s+1]}(\p_r-2p_{\s+1})-\sum_{r\in[1,\s+1]}(2p_r-2p_{\s+1})=\sum_{r\in[1,\s+1]}\ps(r).$$
This equals $(1+\ps(\s))+\ps(\s+1)=(1-1)+0=0=\L'_k-\L_k$ if $\s$ is even and equals 
$1+\ps(\s+1)=1-1=0=\L'_k-\L_k$ if $\s$ is odd. (We use 1.6.)

Finally assume that $k>2p_1$. We have $\p_r\le k,2p_r<k$ for all $r\in[0,1+\s]$ hence 
$\L''_k-\L_k=0-0=0=\L'_k-\L_k$.

We see that for any $k\ge0$ we have $\L''_k-\L_k=\L'_k-\L_k$ and (e) follows.

\subhead 3.7\endsubhead
We prove 0.6 in the case where $G$ is as in 1.3 with $\k=0,Q=0,p\ne2$. Let 
$p_*=(p_1\ge p_2\ge\do\ge p_\s)\in\cp_n$.
Let $B,B'$ be as in 2.4; define $u_w$ in terms of the excellent decomposition 2.2(a) of $w=w_{p_*}\i$ as in 2.4. 
Let $\g_{p_*}$ be the $G$-conjugacy class of $u_w$. Let $N_0=u_w\i-1$. By 2.6, $N_0$ has Jordan blocks of sizes
$2p_1,2p_2,\do,2p_\s$. Hence for any $k\ge0$ we have $\dim N_0^k=\sum_{r\in[1,\s]}\max(2p_r-k,0)$. By 2.4(b) we 
have $(B,u_wBu_w\i)\in\co_w$. Since 
$w\in(C_{p_*})_{min}$, we have $C_{p_*}\dsv\g_{p_*}$. Now let $\g'\in\uuG$ be such that $C_{p_*}\dsv\g'$. Then 
there exists $g\in\g'$ such that $gBg\i=u_wBu_w\i$. Define $V_*,V'_*\in\cf$ by $B_{V_*}=B$, 
$B_{V'_*}=u_wBu_w\i$. We have $a_{V_*,V'_*}=w$ (since $(B,u_wBu_w\i)\in\co_w$) and $gV_*=V'_*$ (since 
$gBg\i=u_wBu_w\i$). By 3.5(a) for any $k\ge0$ we have $\dim N^kV\ge\sum_{r\in[1,\s]}\max(2p_r-k,0)$ hence
$\dim N^kV\ge\dim N_0^kV$. It follows that the conjugacy class of $u_w\i$ in $GL(V)$ is contained in the closure 
of the conjugacy class of $g$ in $GL(V)$. Since $p\ne2$ it follows that the conjugacy class of $u_w\i$ in $G$ is 
contained in the closure of the conjugacy class of $g$ in $G$. Since $\g_{p_*}=\g_{p_*}\i$ we see that $\g_{p_*}$
is contained in the closure of $\g'$. We see that property $\Pi_{C_{p_*}}$ holds with 
$\g_{C_{p_*}}=\g_{p_*}$. The map $p_*\m\g_{p_*}$ is clearly injective. If $\g_0$ is a distinguished unipotent
class in $G$ then all its Jordan blocks have even sizes (with multiplicity one) hence $\g_0$ is of the form 
$\g_{p_*}$ for some $p_*\in\cp_n$. This completes the proof of 0.6 in our case.

\subhead 3.8\endsubhead
We prove 0.6 in the case where $G$ is as in 1.3 with $\k=1,Q\ne0,p\ne2$. Let 
$p_*=(p_1\ge p_2\ge\do\ge p_\s)\in\cp_n$. Let $B,B'$ be as in 2.4; define $u_w$ in terms of the excellent 
decomposition 2.2(a) of $w=w_{p_*}\i$  as in 2.4. Let $\g_{p_*}$ be the $G$-conjugacy class of $u_w$. Let 
$N_0=u_w\i-1$. By 2.11, $N_0$ has Jordan blocks of sizes $2p_1+\ps(1),2p_2+\ps(2),\do,2p_\s+\ps(s)$ (and $1$ if 
$\s$ is even). Here $\ps$ is as in 1.6. Hence for any $k\ge0$ we have $\dim N_0^k=\L''_k$ (notation of 3.6). By 
2.4(b) we have $(B,u_wBu_w\i)\in\co_w$. Since 
$w\in(C_{p_*})_{min}$, we have $C_{p_*}\dsv\g_{p_*}$. Now let $\g'\in\uuG$ be such that $C_{p_*}\dsv\g'$. Then 
there exists $g\in\g'$ such that $gBg\i=u_wBu_w\i$. Define $V_*,V'_*\in\cf$ by $B_{V_*}=B$, 
$B_{V'_*}=u_wBu_w\i$. We have $a_{V_*,V'_*}=w$ (since $(B,u_wBu_w\i)\in\co_w$) and $gV_*=V'_*$ (since 
$gBg\i=u_wBu_w\i$). By 3.6(c),(e) for any $k\ge0$ we have $\dim N^kV\ge\L''_k$ hence
$\dim N^kV\ge\dim N_0^kV$. It follows that the conjugacy class of $u_w\i$ in $GL(V)$ is contained in the closure 
of the conjugacy class of $g$ in $GL(V)$. Since $p\ne2$ it follows that the conjugacy class of $u_w\i$ in $G$ is 
contained in the closure of the conjugacy class of $g$ in $G$. Since $\g_{p_*}=\g_{p_*}\i$ we see that $\g_{p_*}$
is contained in the closure of $\g'$. We see that property $\Pi_{C_{p_*}}$ holds with $\g_{C_{p_*}}=\g_{p_*}$. 
This completes the proof of 0.6(i) in our case.

We prove the injectivity in 0.6(ii). Let $\Pi$ be the set of all sequences $\p_1\ge\p_2\ge\do\ge\p_{\s'}$ of 
integers $\ge1$ such that $\p_1+\p_2+\do+\p_{\s'}=\nn$. We define a map $\ph:\cp_n@>>>\Pi$ by

(i) $(p_1\ge p_2\ge\do\ge p_\s)\m(2p_1+\ps(1)\ge2p_2+\ps(2)\ge\do\ge2p_\s+\ps(\s))$ if $\s$ is odd;

(ii) $(p_1\ge p_2\ge\do\ge p_\s)\m(2p_1+\ps(1)\ge2p_2+\ps(2)\ge\do\ge2p_\s+\ps(\s)\ge1)$ if $\s$ is even.
\nl
It is enough to show that $\ph$ is injective.
We show that $2p_i+\ps(i)\ge2p_{i+1}+\ps(i+1)$ for $i\in[1,\s-1]$. If 
$p_i>p_{i+1}$ then $2p_i\ge2p_{i+1}+2$ and it is enough to show that $\ps(i)-\ps(i+1)\ge-2$; this is clear since 
$\ps(i)\ge-1,-\ps(i+1)\ge-1$. So we can assume that $p_i=p_{i+1}$. If $i$ is even then $\ps(i+1)=0$, $\ps(i)\ge0$
hence $2p_i+\ps(i)\ge2p_{i+1}+\ps(i+1)$. If $i$ is odd then $\ps(i+1)\le0$, $\ps(i)=0$ and again 
$2p_i+\ps(i)\ge2p_{i+1}+\ps(i+1)$. Moreover in case (ii) we have $2p_\s+\ps(\s)\ge1$ (since $2p_\s\ge2$,
$\ps(\s)\ge-1$). We see that $\ph$ is well defined. 
Assume now that $p_*=(p_1\ge p_2\ge\do\ge p_\s)\in\cp_n^+$, $p'_*=(p'_1\ge p'_2\ge\do\ge p'_\s)\in\cp_n^+$ 
satisfy $2p_i+\ps(i)=2p'_i+\ps'(i)$ for $i\in[1,\s]$ (here $\ps'$ is defined in terms of $p'_*$ in the same way 
as $\ps$ is defined in terms of $p_*$). Since $\ps(1)=\ps'(1)=1$ we see that $p_1=p'_1$. Assume now that $i\ge2$ 
and that $p_j=p'_j$ for $j\in[1,i-1]$. From our assumption we have $\ps(i)=\ps'(i)\mod2$. If $i$ is odd then 
$\ps(i),\ps'(i)$ belong to $\{0,1\}$ hence $\ps(i)=\ps'(i)$ and $p_i=p'_i$. If $i$ is even then $\ps(i),\ps'(i)$ 
belong to $\{0,-1\}$ hence $\ps(i)=\ps'(i)$ and $p_i=p'_i$. Thus $p_*=p'_*$. Similarly we see that if
$p_*=(p_1\ge p_2\ge\do\ge p_\s)\in\cp_n-\cp_n^+$, $p'_*=(p'_1\ge p'_2\ge\do\ge p'_\s)\in\cp_n-\cp_n^+$ satisfy 
$2p_i+\ps(i)=2p'_i+\ps'(i)$ for $i\in[1,\s]$ then $p_*=p'_*$. This proves the injectivity statement in 0.6(ii). 
Now let $\g_0$ be a distinguished unipotent class in $G$. Then, if $u_0\in\g_0$, the Jordan blocks of $u_0-1$
have sizes $2x_1+1>2x_2+1>\do>2x_f+1$ where $x_1>x_2>\do>x_f$ are integers $\ge0$ and $f$ is odd. Let 
$p_*=(x_1\ge x_2+1\ge x_3\ge x_4+1\ge\do\ge x_{f-2}\ge x_{f-1}+1\ge x_f)$ if $x_f>0$ and
$p_*=(x_1\ge x_2+1\ge x_3\ge x_4+1\ge\do\ge x_{f-2}\ge x_{f-1}+1)$ if $x_f=0$. Then $\g_{p_*}=\g_0$. This 
completes the proof of 0.6 in our case.

\subhead 3.9\endsubhead
We prove 0.6 in the case where $G$ is as in 1.3 with $\k=0,Q\ne0,p\ne2,\nn\ge8$. Let 
$p_*=(p_1\ge p_2\ge\do\ge p_\s)\in\cp^+_n$. Let $w_{p_*}\in W'$ be as in 1.6. Let $w=w_{p_*}\i$.
By the argument in 2.12 we can find $U_*\in\cf$ and a unipotent element $u\in G$ such that $a_{U_*,uU_*}=w$ and 
such that, setting $N_0=u\i-1$, the Jordan blocks of $N_0$ have sizes 
$2p_1+\ps(1),2p_2+\ps(2),\do,2p_\s+\ps(s)$ ($\ps$ as in 1.6). We can also assume that $U_*\in\cf'$.
Note that for any $k\ge0$ we have $\dim N_0^k=\L''_k$ (notation of 3.6).
Let $\g_{p_*}$ be the $G$-conjugacy class of $u$. Let $B=B_{U_*}$. We have $(B,uBu\i)\in\co_w$. 
Since $w\in(C'_{p_*})_{min}$, we have $C'_{p_*}\dsv\g_{p_*}$. Now let $\g'\in\uuG$ be such that 
$C'_{p_*}\dsv\g'$. Then there exists $g\in\g'$ such that $gBg\i=uBu\i$. We set $U'_*=uU_*u\i$. We
have $a_{U_*,U'_*}=w$ and $gU_*=U'_*$ (since $gBg\i=uBu\i$). By 3.6(c),(e) for any $k\ge0$ we have 
$\dim N^kV\ge\L''_k$ hence $\dim N^kV\ge\dim N_0^kV$. It follows that the conjugacy class of $u\i$ in $GL(V)$ is 
contained in the closure of the conjugacy class of $g$ in $GL(V)$. Since $p\ne2$ it follows that the conjugacy 
class of $u$ in $Is(V)$ is contained in the closure of the conjugacy class of $g$ in $Is(V)$. Hence $\g_{p_*}$ is
contained in the closure of $\g'\cup(h\g'h\i)$ where $h\in Is(V)-G$. Then either $\g_{p_*}$ is contained in the 
closure of $\g'$ or $\g_{p_*}$ is contained in the closure of $h\g'h\i$. In the last case we see that
$h\i\g_{p_*}h$ is contained in the closure of $\g'$; but since $2p_1+\ps(1)$ is odd we have 
$h\i\g_{p_*}h=\g_{p_*}$ hence we have again that $\g_{p_*}$ is contained in the closure of $\g'$. We see that 
property $\Pi_{C'_{p_*}}$ holds with $\g_{C'_{p_*}}=\g_{p_*}$. This completes the proof of 0.6(i) in our case. 
The proof of the injectivity in 0.6(ii) is entirely similar to the proof in 3.8. Now let $\g_0$ be a distinguished
unipotent class in $G$. Then, if $u_0\in\g_0$, the Jordan blocks of $u_0-1$ have sizes
$2x_1+1>2x_2+1>\do>2x_f+1$ where $x_1>x_2>\do>x_f$ are integers $\ge0$ and $f$ is even. Let 
$p_*=(x_1\ge x_2+1\ge x_3\ge x_4+1\ge\do\ge x_{f-1}\ge x_f+1)$. Then $\g_{p_*}=\g_0$. This completes the proof of
0.6 in our case.

\head 4. Basic unipotent classes\endhead
\subhead 4.1\endsubhead
In this section there is no restriction on $p$.
Let $G'$ be a connected reductive group over $\CC$ of the same type as $G$ (with the same root datum as $G$). Then
$\uWW$ for $G$ and $G'$ may be identified. Let $\uuG'$ be the set of unipotent classes of $G'$. There is a well
defined map $\p:\uuG'@>>>\uuG$ given by $\p(\g')=\g$ where $\g'$, $\g$ correspond to the same irreducible 
$\WW$-module under the Springer correspondence. This map is injective, dimension preserving. One can show that it
coincides with the map described in \cite{\SPA, III, 5.2}. Let $\ti\Ph':\uWW_{el}@>>>\uuG'$ be the (injective) 
map
$C\m\g_C$ (as in 0.6, for $G'$ instead of $G$). Let $\ti\Ph=\p\ti\Ph':\uWW_{el}@>>>\uuG$, an injective map. (When
$p$ is not a bad prime for $G$ then $\ti\Ph$ is given by $C\m\g_C$ (as in 0.6); this follows 
from the explicit computation of the map $C\m\g_C$, see below.) A unipotent class of $G$ is said to be {\it basic}
if it is in the image of $\ti\Ph$. Let $\uuG_b$ be the set of basic unipotent classes of $G$. Note that $\p$ 
restricts to a bijection $\uuG'_b@>\si>>\uuG_b$. Let $\Ph:\uWW_{el}@>\si>>\uuG_b$ be the restriction of $\ti\Ph$.

\subhead 4.2\endsubhead
We shall need the following definition. Let $V,(,),Q,n,\k,G$ be as in 1.3 and let 
$p_*=(p_1\ge p_2\ge\do\ge p_\s)$ be in $\cp_n$ (if $(1-\k)Q=0$) and in $\cp_n^+$ (if $(1-\k)Q\ne0$). Let 
$\ps$ be as in 1.6. Let $\g_{p_*}$ be the unipotent class in $G$ such that for some/any $g\in\g_{p_*}$ the Jordan
blocks of $g-1$ have sizes

$2p_1\ge 2p_2\ge\do\ge 2p_\s$ if $\k=0,Q=0$ or if $\k=0,Q\ne0,p=2$,

$2p_1\ge 2p_2\ge\do\ge 2p_\s\ge1$ if $\k=1,Q\ne0,p=2$,

$2p_1+\ps(1)\ge2p_2+\ps(2)\ge\do\ge2p_\s+\ps(\s)$ if $\k=1,Q\ne0,p\ne2,\s=\text{odd}$ or if $\k=0,Q\ne0,p\ne2$,

$2p_1+\ps(1)\ge2p_2+\ps(2)\ge\do\ge2p_\s+\ps(\s)\ge1$ if $\k=1,Q\ne0,p\ne2,\s=\text{even}$,
\nl
and such that (if $p=2$):

(a) for any $i\in[1,\s]$ we have $((g-1)^{2p_i-1}x,x)\ne0$ for some $x\in\ker(g-1)^{2p_i}$.
\nl
Note that in each case the unipotent conjugacy class $\g_{p_*}$ is well defined.

\subhead 4.3\endsubhead
We now describe the bijection $\Ph:\uWW_{el}@>\si>>\uuG_b$ for $G$ almost simple of various types.

If $G$ is of type $A_n$ then $\Ph(C_{cox})$ is the regular unipotent class.

If $V,(,),Q,n,\k,G$ are as in 1.3 and $p_*=(p_1\ge p_2\ge\do\ge p_\s)$ is in $\cp_n$ (if $(1-\k)Q=0$) and in 
$\cp_n^+$ (if $(1-\k)Q\ne0$) then $\Ph(C_{p_*})=\g_{p_*}$ (if $(1-\k)Q=0$) and $\Ph(C'_{p_*})=\g_{p_*}$ (if 
$(1-\k)Q\ne0$). 

When $G$ is of exceptional type we use the notation of \cite{\MI}, \cite{\SPS} for the unipotent classes in $G$;
an element $C\in\uWW_{el}$ is specified by indicating the characteristic polynomial of an element of $C$ acting 
on
the reflection representation of $\WW$, a product of cyclotomic polynomials $\Ph_d$ (an exception is type $F_4$
when there are two choices for $C$ with characteristic polynomial $\Ph_2^2\Ph_6$ in which case we use the 
notation $(\Ph_2^2\Ph_6)'$, $(\Ph_2^2\Ph_6)''$ for what in \cite{\GP, p.407} is denoted by $D_4$, $C_3+A_1$). The
notation $d;C;\g$ means that $C\in\uWW_{el},\g\in\uuG_b,\Ph(C)=\g,d=d_C$ (see 0.2). A symbol $dist$ is added when
$\g$ is distinguished for any $p$; a symbol $dist_p$ is added when $\g$ is distinguished only for the specified 
$p$. The values of $d_C$ are taken from \cite{\GP}.

Type $G_2$.

$2;\Ph_6;G_2\quad\text{dist}$
          
$4;\Ph_3;G_2(a_1)\quad\text{dist}$

$6;\Ph_2^2;\tA_1\quad\text{dist}_3$

\mpb

Type $F_4$.

$4;\Ph_{12};F_4 \quad\text{dist}$

$6; \Ph_8; F_4(a_1)\quad\text{dist}$

$8; \Ph_6^2; F_4(a_2)\quad\text{dist}$

$10; (\Ph_2^2\Ph_6)'; B_3$

$10; (\Ph_2^2\Ph_6)''; C_3$

$12; \Ph_4^2; F_4(a_3)   \quad\text{dist}$

$14; \Ph_2^2\Ph_4; C_3(a_1)    \quad\text{dist}_2$

$16;  \Ph_3^2;\tA_2+A_1   \quad\text{dist}_2$

$24;  \Ph_2^4; A_1+\tA_1$

\mpb

Type $E_6$.

$6;   \Ph_3\Ph_{12}; E_6\quad\text{dist}$

$8;  \Ph_9; E_6(a_1)\quad\text{dist}$

$12; \Ph_3\Ph_6^2; A_5+A_1\quad\text{dist}$

$14; \Ph_2^2\Ph_3\Ph_6; A_5$

$24; \Ph_3^3;2A_2+A_1$

\mpb

Type $E_7$.

$7; \Ph_2\Ph_{18}; E_7\quad\text{dist}$

$9; \Ph_2\Ph_{14}; E_7(a_1)\quad\text{dist}$

$11; \Ph_2\Ph_6\Ph_{12}; E_7(a_2)\quad\text{dist}$

$13; \Ph_2\Ph_6\Ph_{10}; D_6+A_1\quad\text{dist}$

$15;  \Ph_2^3\Ph_{10}; D_6$                 

$17;  \Ph_2\Ph_4\Ph_8; D_6(a_1)+A_1\quad\text{dist}$

$21; \Ph_2\Ph_6^3; D_6(a_2)+A_1\quad\text{dist}$

$23; \Ph_2^3\Ph_6^2; D_6(a_2)$

$25; \Ph_2\Ph_3^2\Ph_6; (A_5+A_1)''$

$31; \Ph_2^5\Ph_6; D_4+A_1$

$33; \Ph_2^3\Ph_4^2; A_3+A_2+A_1$

$63; \Ph_2^7;4A_1$

\mpb

Type $E_8$.

$8; \Ph_{30}; E_8\quad\text{dist}$

$10; \Ph_{24};E_8(a_1)\quad\text{dist}$

$12; \Ph_{20}; E_8(a_2)\quad\text{dist}$

$14; \Ph_6\Ph_{18}; E_7+A_1\quad\text{dist}$

$16; \Ph_{15}; D_8\quad\text{dist}$

$16; \Ph_2^2\Ph_{18}; E_7$

$18; \Ph_2^2\Ph_{14} ; E_7(a_1)+A_1\quad\text{dist}$

$20; \Ph_{12}^2 ; D_8(a_1)\quad\text{dist}$

$22; \Ph_4^2\Ph_{12}; D_7$

$22;  \Ph_6^2\Ph_{12};E_7(a_2)+A_1\quad\text{dist}$

$24; \Ph_{10}^2; A_8\quad\text{dist}$

$24; \Ph_2^2\Ph_6\Ph_{12}; E_7(a_2)$

$26; \Ph_3^2\Ph_{12}; E_6+A_1$

$26;  \Ph_2^2\Ph_6\Ph_{10}; D_7(a_1)\quad\text{dist}_2$

$28; \Ph_3\Ph_9; D_8(a_3)\quad\text{dist}$

$30; \Ph_8^2; A_7\quad\text{dist}_3$

$32; \Ph_2^4\Ph_{10}; D_6$                       

$34; \Ph_2^2\Ph_4\Ph_8; D_5+A_2\quad\text{dist}_2$

$40; \Ph_6^4;2A_4\quad\text{dist}$

$42;  \Ph_2^2\Ph_6^3; A_5A_2$

$44; \Ph_2^4\Ph_6^2; D_6(a_2)$

$44; \Ph_3^2\Ph_6^2; A_5+2A_1$

$46; \Ph_2^2\Ph_3^2\Ph_6;(A_5+A_1)'$

$46; \Ph_2^2\Ph_4^2\Ph_6; D_5(a_1)+A_2$

$48; \Ph_5^2; A_4+A_3$

$60; \Ph_4^4;2A_3$

$64; \Ph_2^6\Ph_6; D_4+A_1$

$66; \Ph_2^4\Ph_4^2; A_3+A_2+A_1$

$80; \Ph_3^4;2A_2+2A_1$

$120; \Ph_2^8;4A_1$ 

\subhead 4.4\endsubhead
We have the following result.

(a) {\it If $\g$ is a distinguished unipotent class of $G$ then $\g$ is a basic unipotent class of $G$.}
\nl
This follows from the known classification of distinguished unipotent classes \cite{\MI}, \cite{\SPA} and the
results in 4.3. For example if $V,(,),Q,n,\k,G$ are as in 1.3 with $p=2$ and $p_*=(p_1\ge p_2\ge\do\ge p_\s)$ is 
in $\cp_n$ (if $(1-\k)Q=0$) and in $\cp_n^+$ (if $(1-\k)Q\ne0$) then $\g_{p_*}$ is distinguished if and only if
for any $j\ge1$ we have $\sha(i\in[1,\s];2p_i=j)\le2$ (and all distinguished classes are of this form).

Next we note the following result:

(b) {\it Let $C\in\uWW_{el}$ and let $g\in\Ph(C)$. If $G$ is semisimple, then $\dim(Z(g))$ is equal to $d_C$ 
(the minimum value of the length function on $C$).}
\nl
When $G$ is almost simple of type $A_n$ this is obvious. When $G$ is almost simple of exceptional type this 
follows from the results in 4.3 and from \cite{\MI,\SPS}. Now assume that $V,(,),Q,n,\k,G$ are as in 1.3 and 
$p_*=(p_1\ge p_2\ge\do\ge p_\s)$ is in $\cp_n$ (if $(1-\k)Q=0$) and in $\cp_n^+$ (if $(1-\k)Q\ne0$). Let 
$d'=\dim Z(g)$. Using 2.2, 2.3 we see that it is enough to show that 

(c) $d'=2\sum_{v=1}^{\s-1}vp_{v+1}+n$ (if $(1-\k)Q=0$) and $d'=2\sum_{v=1}^{\s-1}vp_{v+1}+n-\s$ (if 
$(1-\k)Q\ne0$).
\nl
Since $\dim(\p(\g))=\dim\g$ for any $\g\in\uuG'$ (notation of 4.1) we see that it is enough to prove (c) in the
case where $p=2$. Using the exceptional isogeny from type $B_n$ to type $C_n$ we see that (c) in type $B_n$ 
follows from (c) in type $C_n$. Using \cite{\SPA, II, 6.4, 6.5} we see that (c) in type $D_n$ follows from (c) in 
type $C_n$. Thus we may assume that $\k=0,Q=0,p=2$. For $j\ge1$ let $f_j=\sha(i\in[1,\s];2p_i\ge j)$. By 
\cite{\SPA, II, 6.3, 6.5} we have 
$$d'=\sum_{h\ge1}(f_{2h}^2-f_{2h})+n$$
since $f_{2h-1}=f_{2h}$. It remains to prove the identity $X=2Y$ where 
$$X=\sum_{h\ge1}(f_{2h}^2-f_{2h}),\qua Y=p_2+2p_3+\do+(\s-1)p_\s.$$
We can find integers $a_1,a_2,\do,a_t,b_1,b_2,\do,b_t$ (all $\ge1$) such that 
$p_i=a_1+a_2+\do+a_t$ for $i\in[1,b_1]$, $p_i=a_1+a_2+\do+a_{t-1}$ for $i\in[b_1+1,b_1+b_2]$, $\do$, $p_i=a_1$ for
$i\in[b_1+b_2+\do+b_{t-1}+1,b_1+b_2+\do+b_t]$. We have 
$$\align&X=a_1((b_1+b_2+\do+b_t)^2-(b_1+b_2+\do+b_t))\\&+a_2((b_1+b_2+\do+b_{t-1})^2-(b_1+b_2+\do+b_{t-1}))+\do
+a_t(b_1^2-b_1),\endalign$$
$$\align&Y=(a_1+a_2+\do+a_t)(b_1^2-b_1)/2\\&+(a_1+a_2+\do+a_{t-1})((b_2^2-b_2)/2-(b_1^2-b_1)/2)+\do+\\&
a_1((b_t^2-b_t)/2-(b_{t-1}^2-b_{t-1})/2).\endalign$$
The equality $X=2Y$ follows. This completes the proof of (b).

\subhead 4.5\endsubhead
We now define a map $\Ph:\uWW@>>>\uuG$ extending the map $\Ph:\uWW_{el}@>>>\uuG_b$ in 4.1. Let $C\in\uWW$. We can 
find $J\sub S$ and an elliptic conjugacy class $D$ of the Weyl group $\WW_J$ such that $D=C\cap\WW_J$. Let $P$ 
be a parabolic subgroup of $G$ of type $J$. Let $L$ be a Levi subgroup of $P$. Let $\g_D=\Ph_L(D)$, a unipotent 
class of $L$ (here $\Ph_L$ is the map $\Ph$ of 4.1 with $G,\WW$ replaced by $L,\WW_J$). Let $\g$ be the unipotent
class of $G$ containing $\g_D$. We set $\Ph(C)=\g$. We show that $\g$ is independent of the choices made. Assume 
that we have also $D'=C\cap\WW_{J'}$ where $J'\sub S$ and $D'$ is an elliptic conjugacy class of the Weyl group
$\WW_{J'}$. Let $P'$ be a parabolic subgroup of $G$ of type $J'$. Let $L'$ be a Levi subgroup of $P$. Let 
$\g_{D'}=\Ph_{L'}(D')$, a unipotent class of $L'$. Let $\g'$ be the unipotent class of $G$ containing $\g_{D'}$. 
We must show that $\g=\g'$. By \cite{\GP, 3.2.12} there exists $x\in\WW$ such that $xJx\i=J'$ and $xDx\i=D'$. We 
can find an element $\dx\in G$ such that $\dx L\dx\i=L'$ and such that conjugation by $\dx$ induces the 
isomorphism $\WW_J@>>>\WW_{J'}$ given by $w\m xwx\i$. By functoriality we must have 
$\dx\Ph_L(D)\dx\i=\Ph_{L'}(D')$. It follows that $\dx\g\dx\i=\g'$ hence $\g=\g'$. We see that $C\m \Ph(C)$ is a
well defined map $\uWW@>>>\uuG$; it clearly extends the map $\Ph:\uWW_{el}@>>>\uuG_b$ in 4.1. Also, if $p$ is not
a bad prime for $G$ then the map $\uWW@>>>\uuG$ just defined coincides with the map $\uWW@>>>\uuG$ given by 0.4. 
(This follows from 1.1.)

Note that $\Ph:\uWW@>>>\uuG$ can be described explicitly for any $G$ using the description of the bijections 
$\uWW_{el}@>\si>>\uuG_b$ given in 4.3 (with $G$ replaced by a Levi subgroup of a parabolic subgroup of $G$). 

We show that 

(a) {\it $\Ph:\uWW@>>>\uuG$ is surjective.}
\nl
Let $\g\in\uuG$. We can find a parabolic subgroup $P$ of $G$ with Levi subgroup $L$ and a
distinguished unipotent class $\g_1$ of $L$ such that $\g_1\sub\g$. Let $J$ be the subset $S$ such that $P$ is of
type $J$. By 4.4(a), $\g_1$ is a basic unipotent class of $L$. Hence we can find an elliptic conjugacy class 
$D$ of $\WW_J$ such that $\Ph_L(D)=\g_1$ ($\Ph_L$ is the map $\Ph$ of 4.1 with $G,\WW$ replaced by $L,\WW_J$). Let
$C$ be the conjugacy class in $\WW$ that contains $D$. By the arguments above we have $\Ph(C)=\g$. This proves 
(a).

\subhead 4.6\endsubhead
In this subsection we show, assuming that all simple factors of $G$ are of type $A_n,B_n,C_n$ or $D_n$, that a 
part of Theorem 0.4(i) holds for the map $\Ph:\uWW@>>>\uuG$ even in bad characteristic.

(a) {\it Let $C\in\uWW$ and let $\g=\Ph(C)\in\uuG$. Then $C\dsv\g$.}
\nl
As in 1.1 we can assume that $C\in\uWW_{el}$. If $p\ne2$ the result follows from 0.6. We now assume that $p=2$. We
can also assume that $G$ is almost simple of type $\ne A_n$. We can now assume that $V,(,),Q,n,\k,G$ are as in 1.3
with $p=2$ and that $C=C_{p_*}$ (if $(1-\k)Q=0$) and $C=C'_{p_*}$ (if $(1-\k)Q\ne0$) where 
$p_*=(p_1\ge p_2\ge\do\ge p_\s)$ is in $\cp_n$ (if $(1-\k)Q=0$) and in $\cp_n^+$ (if $(1-\k)Q\ne0$). Using the 
exceptional bijection from type $B_n$ to type $C_n$ we see that the result for type $B_n$ follows from the result
in type $C_n$. Thus we can assume in addition that $\k=0$. As in the proof in 3.7, 3.9 we see that there exists 
$\g\in\uuG$ such that if $g\in\g$ then $g-1$ has Jordan blocks of sizes $2p_1\ge2p_2\ge\do\ge2p_\s$ and $C\dsv\g$.
It remains to show that $g$ satisfies the conditions 4.2(a). It is enough to show that there exists a direct sum 
decomposition $V=V^1\op V^2\op\do V^m$ such that $(V^i,V^j)=0$ for $i\ne j$, $V^i\op(V^i)^\pe=V$ for each $i$ and
such that for each $i$, $V^i$ is $g$-stable and $N:=g-1:V^i@>>>V^i$ has a single Jordan block. To do this we use 
3.5(c) (applied to $V_*\in\cf$ or $V_*\in\cf'$ such that $a_{V_*,gV_*}=w_{p_*}$) and we take for $V^i$ the 
subspaces $X_r$ in 3.4 for $r\in[1,\s]$. This completes the proof of (a).

We expect that (a) holds without assumption on $G$.

\subhead 4.7\endsubhead
Let $C\in\uWW_{el}$ and let $w\in C_{min}$. Define $u_w\in G$ in terms of any excellent decomposition of $w$ as 
in 2.4. We conjecture that $u_w\in\Ph(C)$. This is supported by the computations in Section 2.

\subhead 4.8\endsubhead
Assume that $\kk=\CC$ and consider the bijection $\uuG_b@>\si>>\uWW_{el}$ inverse to $\Ph$. We expect that this
coincides with the restriction of the map $\uuG@>>>\uWW$ defined in \cite{\KL}. (This holds in every case in which
the last map has been computed, see \cite{\SKLC}, \cite{\SKL}; in particular it holds when $G$ is as in 1.3.)

\head 5. $C$-small classes\endhead
\subhead 5.1\endsubhead
In this section we fix $C\in\uWW_{el}$. 

Let $\car$ be the reflection representation of $\WW$.
Then $\det(1-w,\car)\in\NN_{>0}$. If $p>1$ we denote by $\det(1-w,\car)^*$ the part prime to $p$ of 
$\det(1-w,\car)$; if $p=0$ we set $\det(1-w,\car)^*=\det(1-w,\car)$. We have the following result.

\proclaim{Theorem 5.2} The isotropy groups of the $G_{ad}$-action 0.2 on $\fB_w$ are finite abelian of order
dividing $\det(1-w,\car)^*$.
\endproclaim
For the proof we shall need the following result which will be proved in 5.3.

(a) {\it If $w',w''\in C_{min}$ then there exists an isomorphism $\fB_{w'}@>\si>>\fB_{w''}$ commuting with the 
$G_{ad}$-actions and commuting with the first projections $\fB_{w'}@>>>G$, $\fB_{w''}@>>>G$.}
\nl
Let $d$ be the order of $w$. 
Using (a) and a result of Geck and Michel \cite{\GP, 4.3.5} we see that we can assume that $w$ is a "good element"
in the sense of \cite{\GP, 4.3.1}. Let $\b^+$ be the braid monoid attached 
to the Coxeter grop $\WW$. Let $w_1\m\hat w_1$ be the canonical imbedding $\WW@>>>\b^+$, see \cite{\GP, 4.1.1}. 
Let $(\hat w)^d$ be the $d$-th power of $\hat w$ in $\b^+$. Let $w_0$ be the longest element of $\WW$. Since $w$ 
is good there exists $z\in\b^+$ such that $(\hat w)^d=\hat w_0z$ in $\b^+$.
Let $s_1s_2\do s_k$ be a reduced expression of $w$. Let $s'_1s'_2\do s'_f$ be a reduced expression of $w_0$.
We can find a sequence $s''_1,s''_2,\do,s''_h$ in $S$ such that $z=\hs''_1\hs''_2\do\hs''_h$. We have
$$(\hs_1\hs_2\do\hs_k)(\hs_1\hs_2\do\hs_k)\do(\hs_1\hs_2\do\hs_k)=\hs'_1\hs'_2\do\hs'_f\hs''_1\hs''_2\do\hs''_h.$$
(The left hand side contains $kd$ factors $\hs_i$. The right hand side contains $f+h$ factors.) We must have
$kd=f+h$. Moreover from the definition of $\b^+$ there exist $\ss^1,\ss^2,\do,\ss^m$ ($m\ge2$) such that each 
$\ss^r$ is a sequence $\ss^r_1,\ss^r_2,\do,\ss^r_{kd}$ in $S$, $\ss^1$ is the sequence 

$s_1,s_2,\do,s_k,s_1,s_2,\do,s_k,\do,s_1,s_2,s_k$,
\nl
($kd$ terms), $\ss^m$ is the sequence 

$s'_1,s'_2,\do,s'_f,s''_1,s''_2,s''_h$
\nl
and for any $r\in[1,m-1]$ the sequence $\ss^{r+1}$ is obtained from the sequence $\ss^r$ by replacing a string 
$\ss^r_{e+1},\ss^r_{e+2},\do,\ss^r_{e+u}$ of the form $s,t,s,t,\do$ ($u$ terms, $s\ne t$ in $S$, $st$ of order 
$u$ in $\WW$) by the string $t,s,t,s,\do$ ($u$ terms). Now let $(g,B)\in\fB_w$, let 
$\fZ=\{c\in G; cgc\i=g,cBc\i=B\}$ and let $c\in\fZ$.
We define a sequence $B_0,B_1,\do,B_{kd}$ in $\cb$ by the following requirements:
$B_{ik}=g^iBg^{-i}$ for $i\in[0,d]$, $(B_{ik+j-1},B_{ik+j})\in\co_{s_j}$ for $i\in[0,d-1],j\in[1,k]$.
This sequence is uniquely determined. Now conjugation by $c$ preserves each of 
$B,gBg\i,g^2Bg^{-2},\do,g^dBg^{-d}$ hence (by uniqueness) it automatically preserves each $B_v$, $v\in[0,kd]$. 
We define a sequence $B_*^1,B_*^2,\do,B_*^m$ such that each $B_*^r$ is a sequence $(B_0^r,B_1^r,\do,B_{kd}^r)$ in
$\cb$ satisfying $(B^r_{j-1},B^r_j)\in\co_{\ss^r_j}$ for $j\in[1,kd]$, as follows: 
$B_*^1=(B_0,B_1,\do,B_{hd})$ and for $r\in[1,m-1]$, $B_*^{r+1}$ is obtained from $B_*^r$ by replacing the string 
\lb $B^r_e,B^r_{e+1},\do,B^r_{e+u}$ (where 

$(\ss^r_{e+1},\ss^r_{e+2},\do,\ss^r_{e+u})=(s,t,s,t,\do)$
\nl
as above) by the string $B^{r+1}_e,B^{r+1}_{e+1},\do,B^{r+1}_{e+u}$ defined by 
$$\align&B^{r+1}_e=B^r_e, B^{r+1}_{e+u}=B^r_{e+u}, (B^{r+1}_e,B^{r+1}_{e+1})\in\co_t, 
(B^{r+1}_{e+1},B^{r+1}_{e+2})\in\co_s,\\& (B^{r+1}_{e+2},B^{r+1}_{e+3})\in\co_t, \do.\endalign$$
(Note that $B^{r+1}_e,B^{r+1}_{e+1},\do,B^{r+1}_{e+u}$ are uniquely determined since
$(B^r_e,B^r_{e+u})\in\co_{stst\do}=\co_{tsts\do}$ and $stst\do,tsts\do$ are reduced expressions in $\WW$.) We 
note that for any $r\in[1,m]$ any Borel subgroup in the sequence $B^r_*$ is stable under conjugation by $c$. (For 
$r=1$ this has been already observed. The general case follows by induction on $r$ using the uniqueness in the 
previous sentence.) In particular any Borel subgroup in the sequence $B^m_*$ is stable under conjugation by $c$. 
From the definitions we see that $(B^m_0,B^m_f)\in\co_{w_0}$ that is, $B^m_0,B^m_f$ are opposed Borel subgroups. 
Since both are stable under conjugation by $c$ we see that $c$ belongs to $B^m_0\cap B^m_f$, a maximal torus 
independent of $c$. We see that $\fZ$ is contained in the torus $B^m_0\cap B^m_f$. Hence $\fZ$ is a 
diagonalizable group.

Now if $c\in\fZ$ then $cBc\i=B$ and $cgBg\i c\i=gBg\i$ hence $c\in B\cap gBg\i$. Thus $\fZ$ is a diagonalizable
subgroup of the connected solvable group $B\cap gBg\i$. Hence we can find a maximal torus $T$ of $B\cap gBg\i$ 
such that $\fZ\sub T$. We can find $a,a'$ in $U_B$ and $y$ in the normalizer of $T$ such that $g=aya'$; moreover 
$y$ is uniquely determined. For $c\in\fZ$ we have $g=cgc\i=cac\i cyc\i ca'c\i$. Since $c\in T$ we see that
$cac\i,ca'c\i$ belong to $U_B$ and $cyc\i$ belongs to the normalizer of $T$. By the 
uniqueness statement above we see that $cyc\i=y$. Thus $\fZ$ is contained in $T^y$ (the fixed point set of 
$\Ad(y):T@>>>T$). Let $\bar\fZ,\bar T,\bar y$ be the image of $\fZ,T,y$ in $G_{ad}$. Then $\bar\fZ$ is contained
in $\bar T^{\bar y}$ (the fixed point set of $\Ad(\bar y):\bar T@>>>\bar T$). Since the conjugacy class of $w$ is
elliptic and  $(B,gBg\i)\in\co_w$ (that is $(B,yBy\i)\in\co_w$) we see that $\bar T^{\bar y}$ is a finite 
abelian group of order $\det(1-w,\car)^*$. Hence $\bar\fZ$ is a finite abelian group of order dividing 
$\det(1-w,\car)^*$. This completes the proof

\subhead 5.3\endsubhead
We prove 5.2(a). By \cite{\GP, 3.2.7(P2)}, there exists a sequence 

$w'=w_1,w_2,\do,w_n=w''$
\nl
in $W$ such that for $i\in[1,n-1]$ we have $w_i=b_ic_i,w_{i+1}=c_ib_i$ where $b_i,c_i$ in $W$ satisfy 
$\ul(b_i)+\ul(c_i)=\ul(b_ic_i)=\ul(c_ib_i)$. Hence we may assume that $w'=bc,w''=cb$ where $b,c$ in $W$ satisfy 
$\ul(b)+\ul(c)=\ul(bc)=\ul(cb)$. 
If $(g,B)\in\fB_{bc}$ then there is a unique $B'\in\cb$ such that $(B,B')\in\co_b$, $(B',gBg\i)\in\co_c$. We have
$(gBg\i,gB'g\i)\in\co_b$ hence $(B',gB'g\i)\in\co_{cb}$ so that $(g,B')\in\fB_{cb}$. Thus we have defined a 
morphism $\a:\fB_{bc}@>>>\fB_{cb}$, $(g,B)\m(g,B')$. Similarly if $(g,B')\in\fB_{cb}$ there exist a unique 
$B''\in\cb$ such that $(B',B'')\in\co_c,(B'',gB'g\i)\in\co_b$; we have $(g,B'')\in\fB_{bc}$. Thus we have defined
a morphism $\a':\fB_{cb}@>>>\fB_{bc}$, $(g,B')\m(g,B'')$. From the definition it is clear that 
$\a'\a(g,B)=(g,gBg\i)$ for all $(g,B)\in\fB_{bc}$ and $\a\a'(g,B')=(g,gB'g\i)$ for all $(g,B')\in\fB_{cb}$. It 
follows that $\a,\a'$ are isomorphisms. They have the required properties.

\subhead 5.4\endsubhead
We give an alternative proof of Theorem 5.2 in the case where $V,(,),Q,n,\k,G$ are as in 1.3. By 5.2(a) we can 
assume that $w=w_{p^*}$ where $p_*$ is in $\cp_n$ (if $(1-\k)Q=0$) and in $\cp_n^+$ (if $(1-\k)Q\ne0$). Let 
$V_*\in\cf$ (if $(1-\k)Q=0$), $V_*\in\cf'$ (if $(1-\k)Q\ne0$) and let $g\in G$ be such that 
$a_{V_*,V'_*}=w_{p_*}$ where $V'_*=gV_*$. It is enough to show that there exists a finite subgroup $\G$ of $G$ 
such that, if $x\in G$ satisfies $xgx\i=g,xV_*=V_*$ then $x\in\G$. Consider the basis $\b$ of $V$ associated in 
3.3(vi) to $V_*,V'_*,g$. Since $\b$ is canonically defined (up to multiplication by $\pm1$) by $V_*,V'_*,g$ we 
see that the analogous basis associated to $xV_*,xV'_*,xgx\i$ is equal to $x\b$ (up to multiplication by $\pm1$).
Since $(xV_*,xV'_*,xgx\i)=(V_*,V'_*,g)$ we see that $x\b$ is equal to $\b$ (up to multiplication by $\pm1$). Let 
$\G$ be the set of all elements of $G$ which map each element of $\b$ to $\pm1$ times itself (a finite abelian 
$2$-subgroup of $G$). We see that $x\in\G$, as required.

\subhead 5.5\endsubhead
Let $\g$ be a conjugacy class of $G$. Let $w\in C_{min}$.
Since $\fB^\g_w$ is a union of $G_{ad}$-orbits in $\fB_w$ and each of these orbits has dimension
equal to $\dim G_{ad}$ (see 5.2) we see that if $\fB^\g_w\ne\em$ then $\dim\fB^\g_w\ge\dim(G_{ad})$.
We say that $\g$ is {\it $C$-small} if $C\dsv\g$ and $\dim\fB_w^\g=\dim(G_{ad})$. (The last condition is
independent of the choice of $w$ in $C_{min}$.)

For any Borel subgroup $B$ of $G$ let $\Om_B$ be the $B-B$ double coset of $G$ such that $(B,xBx\i)\in\co_w$
for some/any $x\in\Om_B$. Assume that $\Om_B\cap\g\ne\em$. We have the following result:

(a) {\it We have $\dim(\Om_B\cap\g)\ge\dim(B/Z_G)$; moreover equality holds if and only if $\g$ is $C$-small.}
\nl
Indeed, we have a fibration $\fB_w^\g@>>>\cb$, $(g,B')\m B'$ whose fibre at $B'$ is $\Om_{B'}\cap\g$. It follows 
that $\dim(\Om_B\cap\g)=\dim\fB_w^\g-\dim\cb$. It remains to use that $\dim(G_{ad})-\dim\cb=\dim(B/Z_G)$.

\proclaim{Corollary 5.6}Let $w\in C_{min}$. Let $g\in G$. The submanifolds $\{(B,B')\in\cb\T\cb;B'=gBg\i\}$ and 
$\co_w$ of $\cb\T\cb$ intersect transversally.
\endproclaim
In the case where $g$ is regular semisimple this is proved in \cite{\LU, 1.1} (without assumption on $w$). As in 
that proof it is enough to verify the following statement:

(a) {\it if $B\in\cb$ and $(B,gBg\i)\in\co_w$ then $(1-\Ad(g))(\fg)+\fb=\fg$.}
\nl
Here $\fb,\fg$ are the Lie algebras of $B,G$. We can assume that $G=G_{ad}$ and that $w$ is good.  Let $\fg^*$ 
be the dual space of $\fg$. We give two proofs (the second one applies only for $p=0$).

For any subspace $\cv$ of $\fg$ let $\cv^\pe$ be the annihilator of $\cv$ in $\fg^*$.
It is enough to show that $\ker((1-\Ad(g)):\fg^*@>>>\fg^*)\cap\fb^\pe=0$. We argue as in the proof of 5.2. 
Let $\x\in\ker((1-\Ad(g)):\fg^*@>>>\fg^*)\cap\fb^\pe$. Let $d,k,f$, 
$\ss^r_j,B^r_j$ $(r\in[1,m],j\in[0,kd])$, $m\ge2$ be as in the proof of 5.2. Let $\fb_r^j$ be the Lie algebra
of $B_r^j$. Let $\x\in\fa$. Since $\Ad(g)^i\x=\x$ for all $i$ and $\x\in\fb^\pe$ we see that
$\x\in(\Ad(g)^i\fb)^\pe$ for $i\in[0,d]$ hence $\x\in(\fb^1_{ik})^\pe$ for $i\in[0,d]$. 
Using the definition of $B^1_j$ for $j\in[0,kd],j\n k\NN$ we see that
$\x\in(\fb^1_j)^\pe$ for any $j\in[0,kd]$. Using the definitions we see by induction on $r$ that
$\x\in(\fb^r_j)^\pe$ for any $j\in[0,kd],r\in[1,m]$. In particular we have
$\x\in(\fb^m_0)^\pe\cap(\fb^m_f)^\pe$. The last intersection is $0$ since $B^m_0,B^m_f$ are opposed Borel 
subgroups. Thus, $\x=0$, as desired.

In the second proof (with $p=0$) let $\fn$ be the Lie algebra of $U_B$. We identify
$\fg=\fg^*$ and $\fn=\fb^\pe$ using the Killing form; we see that it is enough to show that 
$\ker((1-\Ad(g)):\fg@>>>\fg)\cap\fn=0$. The last intersection 
is the Lie algebra of $Z(g)\cap U$. By 5.2, $Z(g)\cap B$ is a finite group. Hence $Z(g)\cap U$ is a finite 
subgroup of $U$ hence $Z(g)\cap U=\{1\}$ and the desired result follows.

\proclaim{Corollary 5.7} We preserve the setup of 5.6.

(i) The variety $\cb_g^w:=\{B\in\cb;(B,gBg\i)\in\co_w\}$ is smooth of pure dimension $\ul(w)$.

(ii) For any $\g\in\uG$, the variety $\fB^\g_w$ is smooth of pure dimension $\dim\g+\ul(w)$.

(iii) Assume that $\g\in\uG$ and $\fB^\g_w\ne\em$. Then for $g\in\g$ we have $\dim(Z(g)/Z_G)\le\ul(w)$ and
$\dim\g\ge\dim(G_{ad})-\ul(w)$.

(iv) With the assumptions of (iii), $\g$ is $C$-small if and only if $\dim(Z(g)/Z_G)=\ul(w)$ that is, if and only 
if $\dim\g=\dim(G_{ad})-\ul(w)$.
\endproclaim
(We use the convention that the empty variety has dimension $d$ for any $d$.)
The variety in (i) may be identified with the intersection in 5.6 (the submanifolds in 5.6 have pure dimension 
$\dim\cb,\dim\cb+\ul(w)$ and $\cb\T\cb$ has dimension $2\dim\cb$); (i) follows. Now (ii) follows from (i) since 
$\fB^\g_w$ is fibred over $\g$ with fibres as in (i) with $g\in\g$. We prove (iii). By 5.2 every $G_{ad}$-orbit in
$\fB^\g_w$ has dimension equal to $\dim(G_{ad})$. Hence from (ii) we see that $\dim\g+\ul(w)\ge\dim(G_{ad})$ and 
(iii) follows. The proof of (iv) is similar to that of (iii). 

\subhead 5.8\endsubhead
We show:

(a) {\it Assume that $\g\in\uG$ is $C$-small. Then the $G_{ad}$-action on $\fB^\g_w$ has finitely many orbits.
Also, if $g\in\g$, the $Z(g)/Z_G$-action on $\cb_g^w$ (by conjugation) has finitely many orbits.}
\nl
We have $\dim\fB^\g_w=\dim(G_{ad})$ (see 5.7(ii)) and every $G_{ad}$-orbit in $\fB^\g_w$ has dimension equal to 
$\dim(G_{ad})$ (see 5.2); the first statement of (a) follows. If $g\in\g$ then $\dim\cb_g^w=\dim(Z(g)/Z_G)$ (by
our assumption and 5.7(i)). The isotropy groups of the $Z(g)/Z_G$-action on $\cb_g^w$ are finite (by 5.2) hence 
every $Z(g)/Z_G$-orbit in $\cb_g^w$ has dimension $\dim(Z(g)/Z_G)$; the second statement of (a) follows. This 
proves (a).

In the following result we assume that $p$ is not a bad prime for $G$. We give an alternative characterization of
$\Ph(C)$ for $C\in\uWW_{el}$ which does not involve the partial order of $\uuG$.

(b) {\it There is a unique $C$-small unipotent conjugacy class in $G$ namely $\Ph(C)$.}
\nl
The fact that $\Ph(C)$ is $C$-small follows from 4.4(b) and 5.8(iv). 
Assume that $\g'\in\uuG$ is any $C$-small unipotent 
class. By $\Pi_C$ we have $\g\sub\bar\g'$. Moreover we have $\dim\g=\dim\g'$ hence $\g=\g'$. This proves (b).
This also proves Theorem 0.7.

\subhead 5.9\endsubhead
Assume that $\kk,F:G@>>>G,F:\cb@>>>\cb,G^F$ are as in the last paragraph of 1.2. Assume that $w\in C_{min}$. Let 
$X_w=\{B\in\cb;(B,FB)\in\co_w\}$, see \cite{\DL}. The finite group $G^F$ acts on $X_w$ by conjugation.
The following result (not used in this paper) is similar to 5.2.

(a) {\it The isotropy groups of the $G^F$-action on $X_w$ are abelian of order prime to $p$.}
\nl
The proof is almost identical to that of 5.2. We can assume that $w$ is good. Let $d,f,s_1,s_2,\do,s_k$, 
$\ss^1,\ss^2,\do,\ss^m$ ($m\ge2$) be as in the proof of 5.2. Let $B\in X_w$, let 
$\cz=\{c\in G^F;cBc\i=B\}=G^F\cap B$ and let $c\in\cz$. We define a sequence $B_0,B_1,\do,B_{hd}$ in $\cb$ by 
the following requirements: $B_{ik}=F^i(B)$ for $i\in[0,d]$, $(B_{ik+j-1},B_{ik+j})\in\co_{s_j}$ for 
$i\in[0,d-1],j\in[1,k]$. This sequence is uniquely determined. Now conjugation by $c$ preserves each of 
$B,FB,F^2B,\do,F^dB$ hence (by uniqueness) it automatically preserves each $B_v$, $v\in[0,kd]$. Starting with 
this sequence and using $\ss^1,\ss^2,\do,\ss^m$ we define a sequence $B_*^1,B_*^2,\do,B_*^m$ ($B_*^r$ is a 
sequence $(B_0^r,B_1^r,\do,B_{kd}^r)$ in $\cb$) as in the proof of 5.2. As in 5.2, any Borel subgroup in the 
sequence $B^r_*$ is stable under conjugation by $c$. In particular $B^m_0,B^m_f$ contain $c$. From the definitions
we see that $(B^m_0,B^m_f)\in\co_{w_0}$ that is, $B^m_0,B^m_f$ are opposed Borel subgroups. We see that $c$ 
belongs to $B^m_0\cap B^m_f$, a maximal torus independent of $c$. Thus $\cz$ is contained in the torus
$B^m_0\cap B^m_f$. This completes the proof of (a).

\enddocument